\documentclass[11pt]{article}

\usepackage{amsmath,amssymb,amsthm,setspace,graphicx,array,tikz,url,blkarray,wasysym,easyReview}
\usepackage[text={6.5in,8.4in}, top=1.3in]{geometry}
\newtheorem{theorem}{Theorem}[section]

\newtheorem{lemma}[theorem]{Lemma}
\newtheorem{corollary}[theorem]{Corollary}
\newtheorem{conjecture}[theorem]{Conjecture}
\newtheorem{problem}[theorem]{Problem}

\newtheorem{remark}[theorem]{Remark}

\newcommand{\B}{\left(\begin{smallmatrix}0 & 0 \\ 0 & 1\end{smallmatrix}\right)}
\newcommand{\M}{\{\left(\begin{smallmatrix}1 & 0\end{smallmatrix}\right),\left(\begin{smallmatrix}1\\0 \end{smallmatrix}\right)\}}
\newcommand{\Flag}{{\mathcal F}}
\newcommand{\Set}{{\mathcal S}}
\newcommand{\Target}{{\mathcal H}}
\newcommand{\Q}{{\mathcal Q}}
\newcommand{\X}{{\mathcal X}}

\begin{document}
	
\title{
	On the maximum density of a matrix and a transcendental Tur\'an-type density
}
	
\author{
	Raphael Yuster
	\thanks{Department of Mathematics, University of Haifa, Haifa 3498838, Israel. Email: raphy@math.haifa.ac.il\;.}
}
	
\date{}
	
\maketitle
	
\setcounter{page}{1}
	
\begin{abstract}
	We prove that the inducibility of $P_4$ in ordered monotone balanced bipartite graphs is $2/e^2$, establishing the smallest known graph with transcendental Tur\'an-type density. Moreover, the limit object is a binary graphon, so it generates a deterministic model. This is a special case of a more general framework addressed here -- the asymptotic maximum density of a constant matrix over an arbitrary symbol set, in a large, possibly monotone, matrix. We solve all $2 \times 2$ monotone cases (one of which corresponds to the aforementioned $P_4$) and all but one of the $2 \times 2$ unrestricted cases.
	
	While $(h!/h^h)^2$ is a lower bound for the asymptotic maximum density of an $h \times h$ matrix,
	we explicitly construct, for all $h \ge 1$, an $h \times h$ minimizer, i.e., a matrix for which this bound is attained. We also sketch how known results on the inducibility of graphs can be modified to show that, as $h$ grows, almost all $h \times h$ $0/1$ matrices are minimizers.
	
\vspace*{3mm}
\noindent
{\bf Keywords:} maximum density; sub-matrix; Tur\'an-type density

\end{abstract}

\section{Introduction}

Tur\'an-type problems are a major area of research in extremal combinatorics that goes back to the fundamental theorem of Tur\'an \cite{turan-1941}. In its generalized form, a
Tur\'an-type problem in a context of combinatorial objects (such as graphs, digraphs, hypergraphs, matrices, ordered graphs, colored graphs) can be stated by 
providing two finite sets of objects: a {\em forbidden set} $\X$ (which may be empty) and a nonempty {\em target set} $\Target$. We then ask for ${\rm ex}(\X, \Target, n)$, the maximum number of induced\footnote{The non-induced case may be reduced to the induced case by enlarging $\Target$.} copies of elements of $\Target$ in a large object
$L$ of order $n$ that does not contain any sub-object isomorphic to an element of $\X$.
If all elements of $\Target$ are of the same fixed order $h$, then we can ask for the {\em Tur\'an density}: the asymptotic maximum probability, as $n$ goes to infinity, that an induced $h$-order sub-object drawn at random from $L$ is isomorphic to an element of $\Target$.
For example, Mantel's Theorem states that ${\rm ex}(\{K_3\}, \{K_2\}, n) = \lfloor n^2/4 \rfloor$ (so the corresponding Tur\'an density is $\frac{1}{2}$).
As another example, for a graph $H$, ${\rm ex}(\emptyset,\{H\},n)$ is the maximum number of induced copies of $H$
in a graph with $n$ vertices and the corresponding maximum asymptotic density is the {\em inducibility} of $H$, an extensively-studied notion. In this paper we consider Tur\'an-type problems in the context of matrices. Our motivation for studying this setting is threefold (i) extend the state of the art on  Tur\'an-type problems in matrices and ordered graphs (ii) study inducibility in the
context of matrices (iii) show that transcendental Tur\'an-type densities arise in this context. We elaborate on related research in each of these three topics in Subsection \ref{subsec:related}.
We next state our main results, which require formalizing some definitions.

Let $H$ be a square matrix of order $h$ and let $M$ be a square\footnote{While the notions mentioned in this paper can be extended to non-square matrices, the technical challenges are already abundant in square matrices, so we prefer to consider the latter in order to reduce notation overload.} matrix, both over an arbitrary set of symbols, assumed without loss of generality to be $\{0,\ldots,k-1\}$ when there are $k$ symbols used. Let $d(H,M)$ be the number of sub-matrices of $M$ equal to $H$, divided by the total number of $h \times h$ sub-matrices of $M$. We call $d(H,M)$ the {\em density} of $H$ in $M$.
One can equivalently define density in graph-theoretic terminology where $H$ and $M$ are edge-colored ordered balanced complete bipartite graphs and we are asking for the probability that a random order-preserving
mapping from $V(H)$ to $V(M)$ induces a subgraph color-isomorphic to $H$.
Notice that for $0/1$ matrices, the analogy is that $H$ and $M$ are ordered balanced bipartite graphs and we are asking for the probability that a random order-preserving
mapping from $V(H)$ to $V(M)$ induces a copy $H$.

Let $f(H,n)$ be the maximum\footnote{A maximum is attained as we can always assume that every symbol appearing in $M$ is a symbol appearing in $H$.} of $d(H,M)$ taken over all $n \times n$ matrices $M$.
Similarly, for a finite set $\X$ of (not necessarily square) matrices, let $f(\X,H,n)$ be the maximum of $d(H,M)$ taken over all $n \times n$ matrices $M$ that do not contain any member of $\X$ as a submatrix; note that $f(H,n)=f(\emptyset,H,n)$.
The Tur\'an-type density parameter of interest is the {\em asymptotic maximum density} of $H$:
$$
f(H) = \lim_{n \rightarrow \infty} f(H,n)  \qquad f(\X,H) = \lim_{n \rightarrow \infty} f(\X,H,n)\;.
~\footnote{By a standard averaging argument, the sequences $f(H,n)$ and $f(\X,H,n)$ are nonincreasing, so the limits exist.}
$$

Suppose now that $H$ is a {\em monotone} matrix, that is
$H(i,j) \le H(i',j')$ for all $i' \ge i$ and $j' \ge j$.
Monotone matrices are well-studied objects in enumerative combinatorics and in particular $0/1$-monotone matrices are the classical {\em Ferrers matrices}. Given $k$, we can associate a
natural {\em forbidden set} $X_k$ with it, namely all row-vectors and $(u,v)$ and column vectors
$(u,v)^T$ such that $k \ge u > v \ge 0$. So, being monotone over the symbol set $\{0,\ldots,k-1\}$ is the same as being $X_{k-1}$-free. Forcing monotonicity on the large matrices, we define
for a monotone $H$ whose largest symbol is $||H||$
$$
f^*(H)= f(X_{||H||},H)\;.
$$
In particular, for a Ferrers matrix $H$ we have $f^*(H)=f(\M,H)$.

The following conjecture seems plausible. It states that for a monotone $H$, the limit object without monotonicity restrictions can be taken to be monotone as well.
\begin{conjecture}\label{conj:monotone}
	Let $H$ be a monotone matrix, then $f(H)=f^*(H)$.
\end{conjecture}

For a given order $h$, one can apply various symmetries to reduce the number of possible $H$ to consider.
We say that $H$ and $H'$ are {\em density isomorphic} if we can obtain $H'$ from $H$ through a sequence of the following operations: (i) reversing the rows, (ii) reversing the columns, (iii) transpose, (iv) bijectively
renaming the symbols. Clearly, if $H$ and $H'$ are density isomorphic, then $f(H,n)=f(H',n)$. Figure \ref{f:1} lists all possible $H$ of order $2$, up to density isomorphism. As can be seen, there are seven isomorphism types, five of which are (density isomorphic to) monotone matrices.
\begin{figure}[!ht]
\[
\begin{pmatrix}	0 & 0 \\ 0 & 0 	\end{pmatrix}~~
\begin{pmatrix}	0 & 0 \\ 0 & 1 	\end{pmatrix}~~
\begin{pmatrix} 0 & 1 \\ 1 & 0 	\end{pmatrix}~~
\begin{pmatrix}	0 & 0 \\ 1 & 1 	\end{pmatrix}~~
\begin{pmatrix}	0 & 0 \\ 1 & 2 	\end{pmatrix}~~
\begin{pmatrix}	0 & 1 \\ 2 & 0 	\end{pmatrix}~~
\begin{pmatrix}	0 & 1 \\ 2 & 3 	\end{pmatrix}
\]
\caption{The $2 \times 2$ matrices, up to density isomorphism.}\label{f:1}
\end{figure}

\begin{figure}[!ht]
\begin{center}
\begin{tikzpicture}[scale=.5]
	\node (a) at (0,2) 		{$\left(\begin{smallmatrix}	0 & 0 \\ 0 & 0 	\end{smallmatrix}\right)$};
	\node (b) at (-2,0) 	{$\left(\begin{smallmatrix}	0 & 0 \\ 0 & 1 	\end{smallmatrix}\right)$};
	\node (c) at (0,0) 		{$\left(\begin{smallmatrix}	0 & 0 \\ 1 & 1 	\end{smallmatrix}\right)$};
	\node (d) at (2,0) 		{$\left(\begin{smallmatrix}	0 & 1 \\ 1 & 0 	\end{smallmatrix}\right)$};
	\node (e) at (-1,-2)	{$\left(\begin{smallmatrix}	0 & 0 \\ 1 & 2 	\end{smallmatrix}\right)$};
	\node (f) at (2,-2) 	{$\left(\begin{smallmatrix}	0 & 1 \\ 2 & 0 	\end{smallmatrix}\right)$};
	\node (g) at (0.5,-4) 	{$\left(\begin{smallmatrix}	0 & 1 \\ 2 & 3 	\end{smallmatrix}\right)$};
	\draw (b) -- (a);
	\draw (c) -- (a);
	\draw (d) -- (a);
	\draw (e) -- (b);
	\draw (e) -- (c);
	\draw (f) -- (d);
	\draw (g) -- (e);
	\draw (g) -- (f);
\end{tikzpicture}
\end{center}
\vspace*{-5mm}
\caption{Hasse diagram of refinement for $2 \times 2$ matrices.}\label{f:2}
\end{figure}
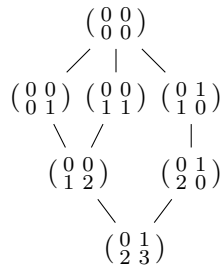

For two matrices $H$ and $H^*$ of the same order, we say that $H^*$ is a {\em refinement} of $H$
if $H(i,j) \neq H(i',j')$ implies $H^*(i,j) \neq H^*(i',j')$ for any two entries $(i,j),(i',j')$.
We observe that refinement is a reflexive partial order and that if $H^*$ is a refinement of $H$, then
$f(H,n) \ge f(H^*,n)$, whence $f(H) \ge f(H^*)$.
Figure \ref{f:2} is the Hasse diagram of refinement for the representatives of the density isomorphism types of $2 \times 2$ matrices.

Let $R=(r_1,\ldots,r_h)$ and $C=(c_1,\ldots,c_h)$ be sequences of nonnegative integers, each summing to $n$.
The {\em $(R,C)$-blowup} of $H$ is the $n \times n$ block matrix $M$ with $h \times h$ blocks, where block $(i,j)$ of $M$ has
$r_i$ rows and $c_i$ columns, and all entries of the block equal $H(i,j)$.
Every sub-matrix of $M$ consisting of a unique entry from each block is equal to $H$, so we have that
$\binom{n}{h}^2 d(H,M) \ge \prod_{i=1}^{h}r_ic_i$.
When each $r_i$ and each $c_i$ are either $\lfloor n/h \rfloor$ or $\lceil n/h \rceil$,
the corresponding $(R,C)$-blowup is called {\em  balanced}\,\footnote{Balanced blowups are not necessarily unique, unless $n$ is a multiple of $h$.}.
If $n$ is a multiple of $h$ and $M$ is the balanced blowup, we have that $\binom{n}{h}^2 d(H,M) \ge (n/h)^{2h}$, thus
\begin{equation}\label{e:minimizer}
f(H) \ge \frac{(h!)^2}{h^{2h}}\;.
\end{equation}
Call $H$ a {\em minimizer} if $f(H)= \frac{(h!)^2}{h^{2h}}$.
Clearly, if $H$ is a minimizer, so are all of its refinements.
The following problem seems intriguing.
\begin{problem}\label{prob:1}
	Given $h$, determine all $h \times h$ minimizers.
\end{problem}
Already for $h=2$, Problem \ref{prob:1} requires work, as can be seen from the proof of our first
result which determines $f^*(H)$ for all monotone $2 \times 2$ matrices and $f(H)$ for all
but one $2 \times 2$ isomorphism type.
\begin{theorem}\label{t:1}
	The values in Table \ref{table:1} hold for $f(H)$ and $f^*(H)$ where $H$ ranges over all isomorphism types of $2 \times 2$ matrices.
\end{theorem}
\begin{table}[ht!]
	\centering
	\renewcommand{\arraystretch}{1.4}
	\begin{tabular}{c||c|c|c|c|c|c|c}
		$H$ & 
		$\left(\begin{smallmatrix}	0 & 0 \\ 0 & 0 	\end{smallmatrix}\right)$ &
		$\left(\begin{smallmatrix}	0 & 0 \\ 0 & 1 	\end{smallmatrix}\right)$ &
		$\left(\begin{smallmatrix}	0 & 0 \\ 1 & 1 	\end{smallmatrix}\right)$ &
		$\left(\begin{smallmatrix}	0 & 1 \\ 1 & 0 	\end{smallmatrix}\right)$ &
		$\left(\begin{smallmatrix}	0 & 0 \\ 1 & 2 	\end{smallmatrix}\right)$ &
		$\left(\begin{smallmatrix}	0 & 1 \\ 2 & 0 	\end{smallmatrix}\right)$ &
		$\left(\begin{smallmatrix}	0 & 1 \\ 2 & 3 	\end{smallmatrix}\right)$ \\
		\hline
		$f(H)$ & $1$ & $\in [\frac{2}{e^2}\,,\,\frac{2}{e^2}+0.001) $  & $\frac{1}{2}$ & $\frac{1}{4}$ & $\frac{1}{4}$ &  $\frac{1}{4}$ & $\frac{1}{4}$\\
		\hline
		$f^*(H)$ & $1$ & $\frac{2}{e^2}$  & $\frac{1}{2}$ & n.a. & $\frac{1}{4}$ &  n.a. & $\frac{1}{4}$
	\end{tabular}
	\caption{The asymptotic maximum density of $2 \times 2$ matrices; n.a. \negthinspace stands for ``non-applicable'' when the matrix is not monotone.}\label{table:1}
\end{table}
An immediate corollary of Theorem \ref{t:1} implies a solution of Problem \ref{prob:1} for $h=2$.
\begin{corollary}\label{coro:1}
	All $2 \times 2$ minimizers are isomorphic to refinements of either $\left(\begin{smallmatrix}0 & 1\\1 & 0\end{smallmatrix}\right)$ or
	$\left(\begin{smallmatrix}0 & 0\\1 & 2\end{smallmatrix}\right)$.
\end{corollary}

\begin{remark}
Theorem \ref{t:1} implies the existence of a {\em transcendental} Tur\'an-type density.
Moreover, all the limit objects that we obtain are {\em $k$-ary} (in particular {\em binary} when the matrix $H$ is $0/1$), i.e., they induce a deterministic model. We note that no explicit transcendental value with an associated deterministic model is yet known for any Tur\'an-type density, assuming a finite forbidden set (we elaborate more on this topic in Subsection \ref{subsec:related}).
\end{remark}
\begin{corollary}
The validity of Conjecture \ref{conj:monotone} implies $f(\B)=2/e^2$\;.
\end{corollary}
Indeed the corollary follows since $f^*(\B)=2/e^2$ by Theorem \ref{t:1}.
Additional support for the conjectured $f(\B)=2/e^2$ is the upper bound for it given in Theorem \ref{t:1}, which is close up to $0.001$. The latter is obtained via flag algebra (see Appendix \ref{sec:flag}).

Returning to Problem \ref{prob:1}, it is not entirely obvious how to prove the existence of an
explicit minimizer for every $h \ge 3$, in particular, proving that the obvious candidate, i.e., the matrix all of whose $h^2$ symbols are distinct (the bottom of the corresponding Hasse diagram) is a minimizer. This is indeed the case, and, as one may expect, far fewer symbols are needed.
\begin{theorem}\label{t:2}
	Let $H$ be the $h \times h$ matrix whose $i$'th row consists entirely of the symbol $i-1$
	except for the last row, whose $j$'th entry is the symbol $h+j-2$. Then $H$ and hence all matrices
	isomorphic to its refinements are minimizers.
\end{theorem}
The number of symbols used in the matrix of Theorem \ref{t:2} is $2h-1$, and it can be made
smaller for all $h \ge 2$ (more details are given in Section \ref{sec:h}; in general, the number of used symbols is $h(1+o_h(1))$).
However, it may be interesting to prove that for all $h \ge 1$, two symbols suffice. While this is trivial for $h=1$ and true for $h=2$ by Theorem \ref{t:1}, we conjecture:
\begin{conjecture}\label{conj:2}
	For all $h \ge 1$, there are $0/1$ $h \times h$ minimizers.
\end{conjecture}
It is not too difficult, using existing techniques, to prove that Conjecture \ref{conj:2} holds for sufficiently large $h$.
In fact, in Section \ref{sec:minimizers}, we sketch how known results about the inducibility of almost
all graphs (see more details in Subsection \ref{subsec:related}) can be used to show that a randomly chosen $0/1$ $h \times h$ matrix is a minimizer with probability going to $1$ as $h$ grows.

An outline of the rest of the paper follows. After a subsection on related research, in Section \ref{sec:exact} we consider all cases of Theorem \ref{t:1} but the case $\B$. The $2/e^2$
proof for $f^*(\B)$, which is also the lower bound for $f(\B)$, is presented in Section \ref{sec:lower}. In Section \ref{sec:h} we prove Theorem \ref{t:2}.
In Section \ref{sec:minimizers} we sketch the proof showing that random $0/1$ matrices are almost surely minimizers.

\subsection{Related research}\label{subsec:related}

Tur\'an-type problems of matrices have been studied since the problem was introduced
by F\"uredi and Hajnal \cite{FH-1992}, who showed that certain extremal problems can be reduced to the study of extremal numbers of non-induced fixed patterns in $0/1$ matrices.
Perhaps the most prominent result in this area is due to Marcus and Tardos \cite{MT-2004} who solved this problem for permutation matrices, consequently establishing the Stanley-Wilf conjecture.
Another notable example is the paper of Pach and Tardos \cite{PT-2006} and the well-known conjecture there, asserting that if the pattern represents the adjacency matrix of a forest, then its Tur\'an number is only slightly larger than linear (it may be $\omega(n \,{\rm polylog} \,n)$ as shown recently by
Pettie and Tardos \cite{PT-2025}).
An interesting {\em induced} Tur\'an-type result in $0/1$ matrices was obtained by Kor\'andi, Pach, and Tomon \cite{KPT-2020} who proved that if $H$ is a $2 \times 2$ matrix which is not homogeneous (not all zero and not all one) and an $n \times n$ matrix $M$ does not contain $H$ as a sub-matrix, then $M$ must contain a $cn \times cn$ homogeneous sub-matrix for some constant $c > 0$.
While Tur\'an-type problems in $0/1$ matrices are equivalent to Tur\'an-type problems in ordered bipartite graphs, we note that the study of extremal problems in ordered (not necessarily bipartite) graphs is also an active
research area (see, e.g., \cite{CFLS-2017,tardos-2018}).

The systematic study of $\operatorname{ex}(\X, \Target, n)$ for $\Target \neq \{K_2\}$ was initiated by
Alon and Shikhelman \cite{AS-2016} although some results in this area date back to papers of Zykov \cite{zykov-1949} and Erd\H{o}s \cite{erdos-1962} determining $\operatorname{ex}(K_t, K_r, n)$ where $t > r$.
The case $\X=\emptyset$, i.e., the study of the maximum induced density for a given graph $H$ in an $n$-vertex graph, and the corresponding maximum asymptotic induced density, a.k.a. {\em inducibility}, denoted here by $i(H)$, was initiated by Pippenger and Golumbic \cite{PG-1975} and has seen a surge of strong results over the last decade. Notice that the definition of $f(H)$ in the present paper
extends the notion of inducibility to matrices, and in particular to ordered bipartite graphs, as the latter correspond to binary matrices. Pippenger and Golumbic observed that
the lexicographic blowup of a given graph $H$ implies that
$i(H) \ge h!/(h^h-h)$ for every graph $H$ on $h$ vertices. We call graphs having their inducibility equal to this generic lower bound {\em weak fractilizers}\footnote{In \cite{FHL-2017}, {\em fractilizers} are graphs for which their nested blowup achieves maximum induced density for every $n$.}.
Notice that weak fractilizers are analogous to our definition of matrix minimizers.
Determining $i(H)$ seems  difficult in general and there are only a handful of graphs for which $i(H)$ is known. While the cases of graphs on three vertices amount to exercises, this is less so for graphs on four vertices \cite{hirst-2014}.
In fact, $i(P_4)$ (here $P_4$ is the the path on four vertices) is still unknown, with record upper and lower bounds appearing in \cite{EL-2015,flagmatic-site}, respectively.
The inducibility of complete bipartite graphs is well understood \cite{BS-1994},
and so is the inducibility of dense Tur\'an graphs \cite{LMZ-2026,LMR-2023,yuster-2026} and the inducibility of graphs that are large blow-ups of other graphs, in terms of their origin \cite{HHN-2014}.
An intriguing conjecture raised in the paper of Pippenger and Golumbic \cite{PG-1975}, which is yet unsolved, is that the cycle $C_h$ for $h \ge 5$ is a weak fractilizer. The case $h=5$ was solved by Balogh, Hu, Lidick\`y, and Pfender \cite{BHLP-2016} with a sophisticated application of flag algebra.
Note that, if true, the Pippenger-Golumbic conjecture would imply that for all $h \ge 5$ there is a weak fractilizer (for $h \le 4$ there is no weak fractilizer, see \cite{EL-2015}; this is unlike the analogous case for $0/1$ minimizers, by Corollary \ref{coro:1}), but at present it is not known
whether for all $h \ge 6$ there is a weak fractilizer, as opposed to our analogous Theorem \ref{t:2} for the matrix setting.
The author proved that almost all graphs $H$ have $i(H)=(1+o_h(1))h!/(h^h-h)$ \cite{yuster-2019} while Fox, Huang, and Lee \cite{FHL-2017} proved that, in fact, almost all graphs $H$ are fractilizers.
Both results use strong asymmetry properties of random graphs, and can be extended to the $0/1$ matrix case, as
we sketch in Section \ref{sec:minimizers}. Inducibility is also studied for directed graphs. In particular, the inducibility of most (but not all) oriented graphs on four vertices is known \cite{BGK-2022,BLPP-2021}.
Finally, a closely related Tur\'an-type inducibility problem for sets, known as the {\em edge-statistics problem} recently attracted attention. Let $\Q_{h,\ell}$ be the set of all graphs on $h$ vertices and $\ell$ edges where $0 < \ell < \binom{h}{2}$. The problem of determining $i(\Q_{h,\ell})$ was raised in \cite{AHKT-2020} who conjectured that $i(\Q_{h,\ell}) \le 1/e+o_k(1)$. Following \cite{KST-2019},
this conjecture was recently solved in \cite{FS-2018,MMNT-2019}.

Many of the recent strong upper bounds (which are sometimes tight) for various Tur\'an-type density problems use the renowned flag algebra method of Razborov \cite{razborov-2007} as their main tool. In particular, flag algebra is used
for computing Tur\'an-type densities of small objects such as graphs, directed graphs, and $3$-uniform hypergraphs \cite{FV-2012,flagmatic-site}. Flag algebra eventually produces a semi-definite optimization program whose coefficients are algebraic combinations of densities in fixed size objects, and does not 
yield transcendental solutions. This suggests that it might not be (even theoretically, regardless of computational power) possible to determine some Tur\'an-type densities using flag algebra and similar ``finite type methods''. While Tur\'an densities can be irrational even if the target set is a single edge \cite{pikhurko-2014}, there are no known transcendental Tur\'an-type densities in any graph-theoretic context\footnote{Pikhurko proved \cite{pikhurko-2014} that if we allow the forbidden set $\X$ to be infinite, then transcendental values must exist. See also the related Question 27 of Fox in \cite{pikhurko-2014}. There are also transcendental densities in the context of fixed patterns in sequences as discovered by Kenyon \cite{kenyon-2026}; the known limit objects of the latter, called {\em sublebesgue measures} induce a random model.}. In fact, Razborov asked in \cite{razborov-2007} Section 6, whether flag algebra
and similar finite type methods suffices to prove all extremal problems involving fixed subgraph densities (see a discussion of Mubayi and Terry in \cite{MT-2019} that exhibited the first extremal graph-theoretic problem whose solution involves explicit transcendental numbers).
Our proof that $f^*(\B) = 2/e^2$ enhances the likelihood that the answer to the aforementioned question is negative, even in a more restrictive setting as we next discuss.

It is well known that every Tur\'an-type graph density corresponds to a limit object \cite{LS-2006}.
For undirected graphs these limit objects are {\em graphons} (analogous limit objects can be defined for edge-colored graphs, for directed graphs, etc.). Formally, a graphon is a symmetric measurable
function $W:[0,1]^2 \rightarrow [0,1]$ and its relation to graph sequences and Tur\'an densities could be seen as follows. To generate an $n$-vertex graph $G_n$ from $W$, we pick $n$ points uniformly at random from $[0,1]$ where the $i'$'th picked point $x_i$ corresponds to vertex labeled $i$ of $G_n$.
We then let $ij \in E(G_n)$ with probability $W(x_i,x_j)$. If $W$ is the graphon corresponding to
a Tur\'an-type density $\alpha$, then with high probability as $n$ goes to infinity, the density of
the target object in the generated graph $G_n$ approaches $\alpha$. Note that, in general, $G_n$ is a random graph even if the $x_i$'s are given (e.g., equally spaced) since $W(x_i,x_j)$ corresponds to a Bernoulli variable. However, if
$W(x,y) \in \{0,1\}$ for all $(x,y) \in [0,1]^2$, the generated graph is deterministically determined
given the selected points. We call such graphons {\em binary} as they generate a deterministic model.
Many (but not all) of the known Tur\'an-type densities computed using flag algebra yield binary graphons.
It seems even more plausible that Tur\'an-type densities which have corresponding binary graphons
could be determined using flag algebra and similar finite type methods. As we shall see however, perhaps surprisingly, our proof for $f^*(\B) = 2/e^2$ yields a corresponding binary graphon and by its transcendental nature, cannot be obtained using flag algebra alone regardless of computational power.

\section{Most $2 \times 2$ matrices}\label{sec:exact}

The ``non-transcendental'' cases of Theorem \ref{t:1} are implied by the following three lemmas.
In all of these lemmas we assume that the containing matrix $M$ is $n \times n$,
it only uses symbols that appear in the considered $H$, that row $i$ of $M$ has $nx_i$ zeros and that $M$ has $xn^2$ zeros. As $xn=\sum_{i=1}^n x_i$, we have
\begin{equation}\label{e:cs}
\sum_{i=1}^n x_i^2+ (1-x_i)^2 \ge n (x^2+(1-x)^2)\;.
\end{equation}
\begin{lemma}\label{l:0110}
	$f(\left(\begin{smallmatrix}0 & 1\\1 & 0\end{smallmatrix}\right))=\frac{1}{4}$\,.
\end{lemma}
\begin{proof}
	We set some notation for various densities of sets of $2 \times 2$ sub-matrices in $M$,
	where the density of a set is the sum of the densities of its elements.
	Let $d_1$ be the density of
	$D_1=\{\left(\begin{smallmatrix}0 & 0\\0 & 0\end{smallmatrix}\right)\,,
	\left(\begin{smallmatrix}1 & 1\\1 & 1\end{smallmatrix}\right)\}$.
	Let $d_2$ be the density of
	$D_2=\{\left(\begin{smallmatrix}0 & 1\\1 & 0\end{smallmatrix}\right)\,,
	\left(\begin{smallmatrix}1 & 0\\0 & 1\end{smallmatrix}\right)\}$.
	Let $d_3$ be the density of
	$D_3=\{\left(\begin{smallmatrix}0 & 1\\0 & 1\end{smallmatrix}\right)\,,
	\left(\begin{smallmatrix}1 & 0\\1 & 0\end{smallmatrix}\right)\}$.
	Let $d_4$ be the density of
	$D_4=\{\left(\begin{smallmatrix}0 & 0\\1 & 1\end{smallmatrix}\right)\,,
	\left(\begin{smallmatrix}1 & 1\\0 & 0\end{smallmatrix}\right)\}$.
	Let $D_5$ be the remaining eight matrices and let the corresponding density be $d_5$. Notice that in each element of $D_5$, precisely three entries are identical.
	We have
	\begin{equation}\label{e:0}
	d_1+d_2+d_3+d_4+d_5=1\;.
	\end{equation}
	Selecting two entries of $M$, not both in the same row and not both in the same column,
	corresponds to selecting a $2 \times 2$ sub-matrix. The probability that both of the selected
	entries are equal is $(1-o_n(1))(x^2+(1-x)^2)$. 
	Since the two diagonals of each $2 \times 2$ matrix have the same chance of being chosen as the
	selected entries, we have that
	\begin{equation}\label{e:1}
	d_1+d_2+\frac{1}{2}d_5 = (x^2+(1-x)^2)(1-o_n(1))\;.
	\end{equation}
	The count of all pairs of entries in $M$ that are both in the same row and are equal is
	$\sum_{i=1}^n (x_i^2+(1-x_i)^2) \binom{n}{2} (1-o_n(1))$.
	We can also count these pairs by considering the densities of all $2 \times 2$ matrices that have at least one
	row in which all pairs are the same. Each pair is counted $(n-1)$ times in this way.
	The matrices in $D_1 \cup D_4$ have two rows with identical entries and the matrices in $D_5$
	have one row with identical entries.
	We therefore have by \eqref{e:cs},
	$$
	\binom{n}{2}^2(2d_1 + 2d_4+d_5) \ge (1-o_n(1))(x^2+(1-x)^2) \frac{n^3}{2}(n-1)
	$$
	and equivalently,
	\begin{equation}\label{e:2}
	d_1 + d_4+\frac{1}{2}d_5  \ge (x^2+(1-x)^2)(1-o_n(1))\;.
	\end{equation}
	Symmetrically, we can consider all pairs of entries in $M$ that are both in the same column and are equal. This yields the inequality
	\begin{equation}\label{e:3}
	d_1 + d_3+\frac{1}{2}d_5  \ge (x^2+(1-x)^2)(1-o_n(1))\;.
	\end{equation}
Maximizing $d_2$ subject to \eqref{e:0}, \eqref{e:1}, \eqref{e:2}, \eqref{e:3} and to $d_i \ge 0$
gives $d_2 \le (1+o_n(1))(1-x^2-(1-x)^2)/2$, so $d_2 \le (1+o_n(1))\frac{1}{4}$
and thus $f(\left(\begin{smallmatrix}0 & 1\\1 & 0\end{smallmatrix}\right)) \le \frac{1}{4}$.
Equality follows from \eqref{e:minimizer}.
\end{proof}
\begin{lemma}\label{l:0012}
	$f(\left(\begin{smallmatrix}0 & 0\\1 & 2\end{smallmatrix}\right))=\frac{1}{4}$\,.
\end{lemma}
\begin{proof}
	Consider two rows $i,j$ of $M$ with $i < j$. We upper-bound the density of $\left(\begin{smallmatrix}0 & 0\\1 & 2\end{smallmatrix}\right)$ in the $2 \times n$ sub-matrix $L$ consisting of these two rows.
	Let $\alpha(1-x_j)$ be the fraction of columns of $L$ equal to $\left(\begin{smallmatrix}0 \\1\end{smallmatrix}\right)$ and let $\beta(1-x_j)$ be the fraction of columns of $L$ equal to $\left(\begin{smallmatrix}0 \\2\end{smallmatrix}\right)$. So, we have that $\alpha+\beta\le 1$,
	$x_i \ge (\alpha + \beta)(1-x_j)$, and the density of $\left(\begin{smallmatrix}0 & 0\\1 & 2\end{smallmatrix}\right)$ in $L$ is
	$$
	2\alpha\beta(1-x_j)^2 \le \frac{(\alpha+\beta)(1-x_j)^2}{2} \le \frac{x_i(1-x_j)}{2}
	$$
	where we have used that $\alpha+\beta \ge 4 \alpha\beta$ since $\alpha+\beta \le 1$. Summing over all pairs $i,j$ with $i < j$ we obtain that
	$$
	f(\left(\begin{smallmatrix}0 & 0\\1 & 2\end{smallmatrix}\right)) \le
	\frac{\binom{n}{2}\sum_{i=1}^n \sum_{j=i+1}^n \frac{1}{2}x_i(1-x_j)}{\binom{n}{2}^2}\;.
	$$ 
	To obtain the lemma, it remains to prove that
	$$
	\sum_{i=1}^n \sum_{j=i+1}^n x_i(1-x_j) \le n^2x(1-x) \le \frac{n^2}{4}\;.
	$$ 
	Indeed,
	$$
	\sum_{i=1}^n \sum_{j=i+1}^n x_i(1-x_j)  = \left(\sum_{i=1}^n x_i(n-i+\frac{x_i}{2})\right) - \frac{n^2x^2}{2} \le \left(n-\frac{xn}{2}\right)xn-\frac{n^2x^2}{2}=n^2x(1-x)\;.
	$$
\end{proof}
\begin{lemma}\label{l:0011}
	$f(\left(\begin{smallmatrix}0 & 0\\1 & 1\end{smallmatrix}\right))=\frac{1}{2}$\,.
\end{lemma}
\begin{proof}
	If $n$ is even, then the balanced blowup has density $\frac{n}{2(n-1)}$, so
	$f(\left(\begin{smallmatrix}0 & 0\\1 & 1\end{smallmatrix}\right)) \ge \frac{1}{2}$.
	Selecting two entries of $M$, not both in the same row and not both in the same column,
	corresponds to selecting a $2 \times 2$ sub-matrix. The probability that one of the selected
	entries is $0$ and the other is not $0$ is $2x(1-x)(1+o_n(1))$, thus $d(\left(\begin{smallmatrix}0 & 0\\1 & 1\end{smallmatrix}\right),M) \le (1+o(1))\frac{1}{2}$ whence
	$f(\left(\begin{smallmatrix}0 & 0\\1 & 1\end{smallmatrix}\right)) \le \frac{1}{2}$.
\end{proof}
The exactly-determined cases of $f(H)$ in Table \ref{table:1} now follow from Lemmas \ref{l:0110}, \ref{l:0012},  \ref{l:0011} and
from the fact that the asymptotic maximum density of a refinement of $H$ is at most the
asymptotic maximum density of $H$. Similarly, all cases of $f^*(H)$ for monotone $H$ other than
$f^*(\B)$ follow immediately since $\frac{1}{4} \le f^*(H) \le f(H)$ and since, as noted in Lemma \ref{l:0011} for even $n$, the balanced blowup of
$\left(\begin{smallmatrix}0 & 0\\1 & 1\end{smallmatrix}\right)$ has density
$\frac{n}{2(n-1)}$ and is monotone.

\section{$\B$}\label{sec:lower}

Our goal in this section is to prove that $f^*(\B) = 2/e^2$.
Since clearly $f(\B) \ge f^*(\B)$, the lower bound for $f(\B)$ in Table \ref{table:1} will follow.
The upper bound $f(\B) < 2/e^2 + 0.001$ is deferred to Appendix \ref{sec:flag}. Thus, together with the results in the previous section, this will conclude the proof of Theorem \ref{t:1}.
Moreover, as stated in the introduction, we will prove that the corresponding limit object
for $f^*(\B)$ is a binary graphon.

let $M$ be a $0/1$ monotone $n \times n$  matrix and let the $i$'th row of
$M$ have $nz_i$ zeros for $1 \le i \le n$. Clearly,
\begin{equation}\label{e:dbm}
d(\B,M)\binom{n}{2}^2 = \sum_{i=1}^n \sum_{j=i+1}^n n^2z_j(z_i-z_j)\;.
\end{equation}
Hence, our goal is to find a nonincreasing sequence $z_1,\ldots,z_n$ where $z_i \le 1$ and $z_n \ge 0$ for which the r.h.s. \negthinspace of the last equality is maximized, as this yields
$f(\M,\B,n)$.
We can view the sequence $z_1,\ldots,z_n$ as a nonincreasing step function with steps of
length $1/n$ in $[0,1)$ in the sense that for the value of the step function at 
$[(i-1)/n,i/n)$ is $z_i$.
In other words, for any given $n$, we are looking for a nonincreasing step function in $[0,1)$ with steps of length $1/n$ and image in $[0,1]$ which maximizes
$$
\sum_{i=1}^n \sum_{j=i+1}^n z_j(z_i-z_j)\;.
$$
Since the integral of any bounded nonincreasing function in $[0,1]$ can be approximated to any desired precision by a nonincreasing step function with steps of length $1/n$, we have that if $g$ is any nonincreasing function in $[0,1]$ with $0 \le g(x) \le 1$,
then
$$
f^*(\B) \ge 4\int_{0}^{1} \int_{x}^{1} g(y)(g(x)-g(y)) \, dy \, dx
$$ 	
where the constant $4$ stems from the term $\binom{n}{2}^2$ in \eqref{e:dbm} which has denominator $4$. Let ${\cal G}^-$ denote the set of all nonincreasing functions in $[0,1]$ with image in $[0,1]$. For $g \in G^-$ consider the functional 
\begin{equation}\label{e:func}
{\mathcal F}(g) = \int_{0}^{1} \int_{x}^{1} g(y)(g(x)-g(y)) \, dy \, dx\;.
\end{equation}
We therefore have that
$$
f^*(\B) = 4\max_{g \in {\mathcal G}^-}{\mathcal F}(g)\;.
$$
It remains to prove the following theorem.
\begin{theorem}\label{t:variational}
	$\max_{g \in {\mathcal G}^-}{\mathcal F}(g) = 1/2e^2$. Consequently, $f^*(\B)=2/e^2$.
	Furthermore, the limit object corresponding to $f^*(\B)$ is a binary graphon.
\end{theorem}
\begin{proof}
	For the lower bound, consider the following function:
\begin{equation}\label{e:g*}
	g^*(x) =
	\begin{cases}
		1 				& {\rm for~} 0 \le x \le \frac{1}{e}\,, \\
		\frac{1}{ex}	& {\rm for~} \frac{1}{e} \le x \le 1\,.
	\end{cases}	
\end{equation}
Clearly $g^* \in {\mathcal G}^-$; in fact $g^*$ is also continuous. We have by \eqref{e:func}
$$
{\mathcal F}(g^*) = \int_{0}^{\frac{1}{e}} \int_{\frac{1}{e}}^{1} \frac{1}{ye}\left(1-\frac{1}{ye}\right) \, dy \, dx
+ \int_{\frac{1}{e}}^{1} \int_{x}^{1} \frac{1}{ye}\left(\frac{1}{xe}-\frac{1}{ye}\right) \, dy \, dx
= \frac{1}{2e^2}\;.
$$
Thus, $\max_{g \in {\mathcal G}^-}{\mathcal F}(g) \ge 1/2e^2$.

We now turn to the upper bound, where we need to show that for every $g \in {\mathcal G}^-$,
${\mathcal F}(g) \le 1/2e^2$. Notice that if we were allowed to assume that $g$ is smooth, then 
this would amount to a smooth variational calculus problem. However, we cannot assume this;
in fact, we cannot even assume that $g$ is continuous as it may potentially have jumps.

Let us first switch to a slightly more convenient form of nondecreasing functions by replacing $g(x)$ with $1-g(x)$ so that that the functional is now
\begin{equation}\label{e:func1}
	{\mathcal F}(g) = \int_{0}^{1} \int_{x}^{1} (1-g(y))(g(y)-g(x)) \, dy \, dx
\end{equation}
and we are maximizing over ${\mathcal G}^+$, the nondecreasing functions in $[0,1]$ so that $g(0) \ge 0$ and $g(1) \le 1$. In fact, we may and will assume that $g(0)=0$ and $g(1)=1$ since this keeps $g$ nondecreasing and changing a single point does not change the integral.

Let us next define some (generalized) functions in order to perform a change of variables.
Let $g^*(x)$ be the Sobolev generalized derivative of $g$. Namely, $g^*$ is such that
$\int_{0}^{x}g^*(y)dy=g(x)$. Notice that $g^*$ is possibly a generalized function.
Indeed, suppose that $g$ has a jump of size $c$ at point $t$. Then $g^*(t)=c\delta_t$ where $\delta_t$ is the usual Dirac delta at point $t$. Further notice that $g^*(x) \ge 0$ as $g$ is nondecreasing. Also observe that $\int_{0}^{1}g^*(x)dx=1$ since $g(1)=1$ (possibly $g^*(1)=c \delta_1$
for some $c > 0$). In other words, $g^*$ is a probability measure on $[0,1]$.
The points where $g^*$ is a Dirac delta correspond to points where the distribution has positive mass.

Let $h=g^{-1}$ and notice that $h$ is a well-defined function by taking $h(u)$ to be the the infimum of
all $x$ for which $g(x) \ge u$ ($h$ is usually referred to in the literature as the quantile function of $g$). Just like $g$, we have that $h$ is nondecreasing in $[0,1]$ and satisfies $h(0)=0$, $h(1)=1$. Thus, its generalized derivative $h^*$ may again be a generalized function, since $h$ may have jumps and just like $g^*$, we have that $h^*$ is a probability measure on $[0,1]$.
Let us now make a change of variables $u=g(x)$ and
$v=g(y)$ so we obtain from \eqref{e:func1}
\begin{equation}\label{e:func2}
	{\mathcal F}(g) = \int_{0}^{1} \int_{u}^{1} (1-v)(v-u)h^*(u)h^*(v) \, dv \, du\;.
\end{equation}
Notice that \eqref{e:func2} is the expected value of $(1-v)(v-u)$ obtained
from two samples from $h^*$ where $u$ is the smaller and $v$ is the larger
and our problem is to find the distribution $h^*$ in $[0,1]$ that maximizes this expected value.
Symmetrically (reflecting the intervals), this is precisely the same distribution which maximizes $u(v-u)$, hence we obtain
\begin{equation}\label{e:func3}
	\max_{g \in {\mathcal G}^+}{\mathcal F}(g) = \max_{h^*} \int_{0}^{1} \int_{u}^{1} u(v-u)h^*(u)h^*(v) \, dv \, du
\end{equation}
where $h^*$ is taken over all probability measures in $[0,1]$.
The maximizing distribution of \eqref{e:func3} has been recently determined by Kenyon \cite{kenyon-2026} using variational calculus, where he used it to determine the
sublebesgue measure of densities of some fixed $0/1$-patterns in large sequences, in particular, the pattern $0101$ which corresponds to \eqref{e:func3}. It turns out that these measures also involve transcendental constants (in such cases, the sublebesgue measure are surprisingly involved, which is unlike our case in which the limit object is a graphon which corresponds to a deterministic model).
An optimal solution for \eqref{e:func2} (see Theorem 4.2 in \cite{kenyon-2026}) is given by:
\begin{equation}\label{e:h*}
	h^*(x) =
\begin{cases}
	0 				& {\rm for~} 0 \le x < \frac{1}{e}\,, \\
	\frac{1}{ex^2}	& {\rm for~} \frac{1}{e} \le x < 1\,, \\
	\frac{\delta_1}{e} & {\rm for~}  x = 1 
\end{cases}	
\end{equation}
for which the integral in \eqref{e:func3} yields $1/2e^2$ as required.
In fact, solving back for $h$ and then for $g$ gives that the corresponding optimizing $g$ is precisely
the function given in \eqref{e:g*} (almost everywhere, i.e., up to choosing $g^*(1)=0$ which, recall, has no effect). For completeness, in Section \ref{sec:var} we present a slightly more verbose description of Kenyon's clever proof of \eqref{e:func3}.

Finally, $g^*$ in \eqref{e:g*} immediately corresponds to our
limit object $W: [0,1]^2 \rightarrow \{0,1\}$ which is defined by $W(x,y)=0$ if $(x,y)$ is below the graph of $g^*$ and $1$ otherwise. Notice that $W$, depicted in Figure \ref{f:graphon}, is indeed a binary graphon, as required.
We recall that to obtain a sequence of monotone binary matrices in which the density of $\B$ converges to $2/e^2$ we may sample the graphon in equally spaced $n$ points $i/n$ for $i=1,\ldots n$
and define the $n \times n$ matrix $A_n$ so that $A_n(i,j) = W(i/n, j/n)$.
\begin{figure}[!ht]
	\includegraphics[scale=0.5,trim=70 140 200 50, clip]{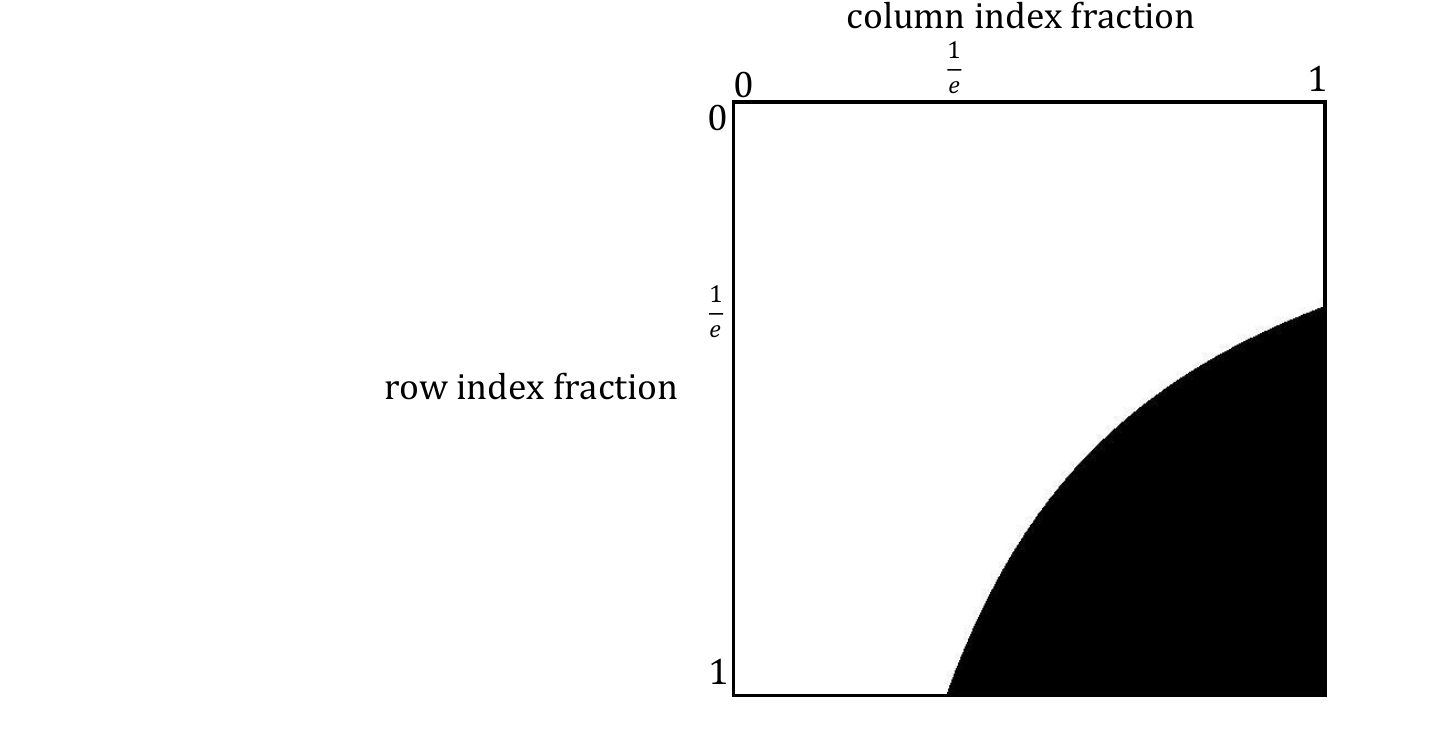}
	\caption{The limit object of a sequence of matrices whose density of $\B$ approaches $2/e^2$. To comply with matrix notation, the vertical axis corresponds to the row index fraction from top to bottom and the horizontal axis corresponds to column index fraction from left to right. The shaded area represents the symbol $1$ and the non-shaded area represents the symbol $0$.}
	\label{f:graphon}
\end{figure} 
\end{proof}

\section{Explicit minimizers for every order}\label{sec:h}

In the proof of Theorem \ref{t:2} we need the following observation.
\begin{lemma}\label{l:inequality}
Suppose that $\alpha_1,\ldots,\alpha_h$ are nonnegative reals with $\alpha_1+\cdots+\alpha_h \le 1$. Then,
$$
(\sum_{i=1}^h \alpha_i)^{h-1} \ge h^h \prod_{i=1}^h \alpha_i\;.
$$
\end{lemma}
\begin{proof}
	Using the AM-GM  inequality, we can substitute $(\sum_{i=1}^h \alpha_i)/h$ with the (not larger) value $(\prod_{i=1}^h \alpha_i)^{1/h}$ and obtain that the lemma's statement reduces to the claim that
	$\prod_{i=1}^h \alpha_i \le 1/h^h$ which clearly holds since $\alpha_1+\cdots+\alpha_h \le 1$.
\end{proof}

\begin{proof}[Proof of Theorem \ref{t:2}]
Recall that $H$ is the $h \times h$ matrix whose $i$'th row consists entirely of the symbol $i-1$
except for the last row, whose $j$'th entry is the symbol $h+j-2$.
We assume that the containing matrix $M$ is $n \times n$ and the set of symbols used in it is $\{0,\ldots,2h-2\}$.
For $0 \le \ell \le h-2$, let $nx_{\ell,i}$ denote the number of occurrences of symbol $\ell$ in row $i$,
and let $x_\ell n^2$ denote the overall number of occurrences of $\ell$ in $M$.
Let $nx_{h-1,i}$ denote the number of occurrences of any of the symbols $h-1,\ldots,2h-2$ in row $i$,
and let $x_{h-1}n^2$ denotes the overall number of occurrences of any of the symbols
$h-1,\ldots,2h-2$ in $M$.
Observe that $\sum_{\ell=0}^{h-1}x_\ell=1$ and that for any row $i$ it holds that $\sum_{\ell=0}^{h-1}x_{\ell,i}=1$. 
We will also assume that $x_{\ell}n$ is an integer
for $0 \le \ell \le h-1$ as this assumption does not change the asymptotic claim of the theorem, since $f(H)$ is defined as a limit where $n$ goes to infinity.

Consider $h$ rows $i_1,\ldots,i_h$ of $M$ with $i_1 < i_2 < \cdots < i_h$.
We upper-bound the density of $H$ in the $h \times n$ sub-matrix $L$ consisting of these rows.

Let $C_j$ be the $j$'th column of $H$. Notice that the bottom entry of $C_j$ is $h+j-2$
and that for $1 \le i \le h-1$, the $i$'th entry of $C_j$ is $i-1$.
Let $\alpha_jx_{h-1,i_h}$ be the fraction of columns of $L$ equal to $C_j$.
So, we have that $\alpha_1+\cdots+\alpha_h \le 1$.
Also observe that 
$x_{\ell-1,i_\ell} \ge (\alpha_1+\cdots+\alpha_h)x_{h-1,i_h}$
for $1 \le \ell \le h-1$. To obtain a copy of $H$ in $L$ we need, in particular, to select precisely
one column of $L$ that equals $C_j$ for every $1 \le j \le h$ (note that such a selection in a necessary condition, but it is not sufficient, as the chosen columns may produce a matrix with the columns of $H$ permuted). Thus, the density of $H$ in $L$ is at most
$$
h! \left(\prod_{i=1}^h \alpha_i \right)x_{h-1,i_h}^h 
\le \frac{h!}{h^h} (\sum_{i=1}^h \alpha_i)^{h-1}x_{h-1,i_h}^h
\le \frac{h!}{h^h} \left(\prod_{\ell=1}^{h-1} x_{\ell-1,i_\ell} \right)x_{h-1,i_h}
\le \frac{h!}{h^h} \left(\prod_{\ell=1}^{h} x_{\ell-1,i_\ell} \right)
$$
where we have used Lemma \ref{l:inequality}.

Summing over all $h$-tuples $i_1,\ldots,i_h$ with $i_1 < i_2 < \cdots < i_h$ we obtain that
$$
d(H,M) \le
\frac{\binom{n}{h}\sum_{i_1=1}^n \sum_{i_2=i_1+1}^n \cdots \sum_{i_h=i_{h-1}+1}^n  \frac{h!}{h^h} \left(\prod_{\ell=1}^{h} x_{\ell-1,i_\ell} \right)}{\binom{n}{h}^2}\;.
$$
To obtain the lemma, it remains to prove that
\begin{equation}\label{e:multisum}
\sum_{i_1=1}^n \sum_{i_2=i_1+1}^n \cdots \sum_{i_h=i_{h-1}+1}^n  \left(\prod_{\ell=1}^{h} x_{\ell-1,i_\ell} \right) \le \frac{n^h}{h^h}\cdot(1+o_n(1))\;.
\end{equation}
Observe that the summation indices of every term in the left hand side of \eqref{e:multisum} correspond to points of the strict upper triangular $h$-dimensional grid $[n]^h$ and that smaller symbols are
associated with lower dimension in the product $\prod_{\ell=1}^{h} x_{\ell-1,i_\ell}$.
We claim that the maximum, up to an additive difference of at most $O(n^{h-1})$, of the left hand side of \eqref{e:multisum} is obtained
when the structure of $M$ is the following {\em canonical form}. All $x_0n$
upper rows entirely consist of the symbol $0$, the next $x_1n$ rows entirely consist of symbol
$1$ and so on, until symbol $h-2$ which occupies $x_{h-2}n$ rows
(notice that this means that the bottom $x_{h-1}n$ rows contain the remaining symbols $h-1,\ldots,2h-2$).
We next prove this claim. In what follows, when comparing symbols, all symbols $h-1,\ldots,2h-2$ are thought of as ``identical'' symbols as if they were all symbol $h-1$ (i.e. all are the ``largest'' symbol).

We will first prove our claim under the assumption that $x_{\ell,i} \in \{0,1\}$ for all variables.
In other words, in every row, all symbols are identical. Notice that this is at least a feasible situation
since we assume that $x_{\ell}n$ is an integer for $0 \le \ell \le h-1$.
So, assume this is the case but that $M$ is not of the canonical form.
Then there must be two consecutive rows, say $m$ and $m+1$, such that $x_{\ell^*,m}=1$
and $x_{\ell,m+1}=1$ yet $\ell < \ell^*$.
We will flip these two rows, claiming that after the flip, the left hand side of \eqref{e:multisum}
does not decrease. Observe that when we do the flip, $x_{\ell,m}$ and $x_{\ell^*,m+1}$ turn from zero to one,
while $x_{\ell,m+1}$ and $x_{\ell^*,m}$ turn from one to zero.
Now, consider the product terms in \eqref{e:multisum} (i.e. terms
of the form  $\prod_{\ell=1}^{h} x_{\ell-1,i_\ell}$) and partition them into
six groups as follows:
(i) Terms that do not involve any of the four variables $x_{\ell,m},x_{\ell^*,m+1},x_{\ell,m+1}, x_{\ell^*,m}$.
(ii) Terms that involve $x_{\ell,m+1}$ and do not involve any of the other three variables.
(iii) Terms that involve $x_{\ell,m}$ and none of the other three variables.
(iv) Terms that involve $x_{\ell^*,m+1}$ and none of the other three variables.
(v) Terms that involve $x_{\ell^*,m}$ and none of the other three variables.
(vi) Terms that involve both $x_{\ell,m}$ and $x_{\ell^*,m+1}$.
Notice that this is indeed a partition, as no term contains two variables associated with the same symbol,
no term contains two variables associated with the same row,
and no term contains both $x_{\ell,m+1}$ and $x_{\ell^*,m}$ (since $\ell < \ell^*$ and $m < m+1$).

Now, following the flip, terms of group (i) do not change their value.
Each term of group (ii) that was equal to $1$ before the flip, corresponds to a unique term of group (iii)
that equals $1$ after the flip.
This is the term obtained by replacing $x_{\ell,m+1}$ with $x_{\ell,m}$. Indeed, in a term of group (ii)
that equals $1$ before the flip, no variable associated with row $m$ exists, since any such variable
must be of the form $x_{t,m}$ where $t < \ell$, but before the flip, $x_{t,m}=0$ since $x_{\ell^*,m}=1$.
Similarly, each term of group (v) that was equal to $1$ before the flip, corresponds to a unique term of group (iv)
that equals $1$ after the flip.
This is the term obtained by replacing $x_{\ell^*,m}$ with $x_{\ell^*,m+1}$. Indeed, in a term of group (v)
that equals $1$ before the flip, no variable associated with row $m+1$ exists, since any such variable
must be of the form $x_{t,m+1}$ where $t > \ell^*$, but before the flip, $x_{t,m+1}=0$ since $x_{\ell,m+1}=1$.
Finally, terms of group (iv) were all zero before the flip. We have proved that the sum of terms
never decreases after performing the flip. Now, repeatedly performing flips we end up with $M$
being of canonical form.

It now remains to prove that given $M$, we can modify it slightly to obtain a matrix such that
in every row, all symbols are the same. Indeed, we will show that this can be done at the expense of
a possible decrease of the sum of the terms of \eqref{e:multisum} that is not larger than $O(n^{h-1})$, which,
in turn, is $o(n^h/h^h)$.
We proceed as follows. Considering the first $h$ rows, take the majority symbol, say $\ell$, in these rows
(there are at least $n$ symbols equal to this $\ell$, by pigeonhole) and make the first row entirely containing the majority symbol $\ell$. Do this by flipping each entry in the first row which is not $\ell$
with an entry which is $\ell$.
Notice that this involves flipping at most $n$ entries, where the distance between the rows of two
flipped symbols is at most $h$. Now, assume that we have already modified $M$ such that the first $m$ rows of $M$ all have the property that in each of these rows, all symbols are the same (distinct rows may have distinct symbols, though). To construct row $m+1$, consider
the majority symbol in rows $m+1,\ldots,m+h$, and use flips to make row $m+1$ entirely contain
that symbol. Once the process ends, the final $M$ has the property that in each row, all symbols are
the same. Also observe that what we have done is a sequence of at most $n^2$ flips, where in each flip,
the two flipped entries are in rows that are at most $h$ apart. It remains to show that each such flip
changes the sum of the terms in \eqref{e:multisum} by at most $O(n^{h-3})$.

Consider such a flip, say of an entry $\ell^*$ in row $m$ with an entry $\ell$ in row $m+t$ where $t < h$.
The variables $x_{\ell,m}$ and $x_{\ell^*,m+t}$ increase by $1/n$ while the variables
$x_{\ell,m+t}$ and $x_{\ell^*,m}$ decrease by $1/n$. As above, partition the terms into six groups.
The definition of groups (i)-(v) is the same except that we now have $m+t$ instead of $m+1$ in the definition of the groups. For group (vi) the definition depends on whether $\ell < \ell^*$ or vice versa. If $\ell < \ell^*$, then the group contains all terms that
involve both $x_{\ell,m}$ and $x_{\ell^*,m+t}$ while if 
$\ell > \ell^*$ then the definition changes to contain all terms that involve both $x_{\ell^*,m}$ and $x_{\ell,m+t}$.
Now, following the flip, terms of group (i) do not change their value.
Terms of group (iii) can only increase their value and so do terms of group (iv).

Consider now terms of group (ii), i.e. those involving $x_{\ell,m+t}$ and do not involve any of the other three variables. Consider such a term $W$. If $\ell = 0$, then $W$ corresponds uniquely to the term
$W'$ where we replace the variable $x_{0,m+t}$ in $W$ with the variable $x_{0,m}$.
Notice that $W'$ is in group (iii) so $W+W'$ is the same before and after the flip.
If $\ell > 0$, then consider the variable $x_{\ell-1,s}$ in $W$. Notice that $s < m+t$. If $s < m$, then
$W$ corresponds uniquely to the term
$W'$ where we replace the variable $x_{\ell,m+t}$ in $W$ with the variable $x_{\ell,m}$.
We again have that $W+W'$ is the same before and after the flip. However, if $m \le s < m+t$
there is no such correspondence. But the number of terms that contain both $x_{\ell,m+t}$ and $x_{\ell-1,s}$
is at most $tn^{h-2}$ (as there are $t$ choices for $s$). As each such term decreases by at most $1/n$,
the net loss for terms of type (ii) is $O(n^{h-3})$.

Analogously,
consider terms of group (v), i.e. those involving $x_{\ell^*,m}$ and do not involve any of the other three variables. Consider such a term $W$. If $\ell^* = h-1$, then $W$ corresponds uniquely to the term
$W'$ where we replace the variable $x_{h-1,m}$ in $W$ with the variable $x_{h-1,m+t}$.
Notice that $W'$ is in group (iv) so $W+W'$ is the same before and after the flip.
If $\ell^* < h-1$, then consider the variable $x_{\ell^*+1,s}$ in $W$. Notice that $s > m$. If $s > m+t$ then
$W$ corresponds uniquely to the term
$W'$ where we replace the variable $x_{\ell^*,m}$ in $W$ with the variable $x_{\ell^*,m+t}$.
We again have that $W+W'$ is the same before and after the flip. However, if $m < s \le m+t$
there is no such correspondence. But the number of terms that contain both $x_{\ell^*,m}$ and $x_{\ell^*+1,s}$
is at most $tn^{h-2}$ (as there are $t$ choices for $s$). As each such term decreases by at most $1/n$,
the net loss for terms of type (v) is $O(n^{h-3})$.

As for terms of type (vi), if $\ell < \ell^*$, then it consists of all terms involving both $x_{\ell,m}$ and $x_{\ell^*,m+t}$, but such terms can only increase after the flip. If $\ell^* < \ell$, then it consists of all terms involving both $x_{\ell^*,m}$ and $x_{\ell,m+t}$. Each such term may decrease by $O(1/n)$, but
there are only $O(n^{h-2})$ such terms, so the net loss for terms of type (vi) is $O(n^{h-3})$.

Finally, for $M$ of canonical form we clearly have
$$
\sum_{i_1=1}^n \sum_{i_2=i_1+1}^n \cdots \sum_{i_h=i_{h-1}+1}^n  \left(\prod_{\ell=1}^{h} x_{\ell-1,i_\ell} \right) \le n^h x_0x_1\cdots x_{h-1} \le \frac{n^h}{h^h}\;.
$$
As we have proved that the canonical form is far from the maximum of the left hand side of \eqref{e:multisum}
by at most $O(n^{h-1})$, we have that \eqref{e:multisum} holds and the theorem follows.
\end{proof}

As noted in the introduction, we can significantly reduce the number of symbols used in an explicitly constructed $h \times h$ minimizer. Recall that $H$ of Theorem \ref{t:2} uses $2h-1$ symbols, but notice that in its bottom row, each symbol is unique in the entire matrix $H$. We only needed this in the proof of Theorem \ref{t:2} to guarantee the following property, thinking of a row of $M$ as a string of length $n$:
The fraction of substrings of length $h$ that equal the bottom row of $H$ is at most $(1+o_n(1))h!/h^h$.
Now, since all $h$ symbols in the bottom row of $H$ are distinct, this clearly holds.
But it is not difficult to prove that we can obtain this property using fewer than $h$ symbols.
To formalize, let $g(h)$ be the least integer $k$ such that there exists a string $w$ of length $h$ over
an alphabet of size $k$ such that for every string of length $n$ it holds that
the fraction of substrings of length $h$ equal to $w$ is at most $(1+o_n(1))h!/h^h$.
By this definition, we can construct, precisely with the same proof of Theorem \ref{t:2}, explicit minimizers using only $h-1+g(h)$ symbols.

It is not difficult to prove that $g(h) \ll h$, but we shall demonstrate this already for the first
nontrivial case $h=3$ where we show that $g(3)=2$ (trivially $g(2)= 2 \le g(3) \le 3$).
Suppose our alphabet is $0/1$. We claim that $w=010$ has the desired property. Consider some
binary string $v=v_1\cdots v_n$. Suppose that $v$ contains $\alpha n$ zeros and $(1-\alpha)n$ ones.
The probability that a random substring of length $3$ of $v$ contains precisely two zeros is at most $3\alpha^2(1-\alpha)$.
Suppose that $v_i=1$ and let $a$ be the number of zeros before $v_i$ in $v$ and $b$ be the number of zeros after $v_i$ in $v$.
Then the number of substrings of $v$ equal to $w=010$ where the $1$ is $v_i$ is $ab$ while the number of $001$ plus the number of $100$
where the $1$ is $v_i$ is $\binom{a}{2}+\binom{b}{2}$. But $\binom{a}{2}+\binom{b}{2} \ge ab - O(n)$ meaning that at most $\frac{1}{2}+o_n(1)$ of the substrings of $v$ of length $3$ with precisely two zeros 
are, in fact, $010$. So, the overall fraction of $010$ substrings is at most $1.5\alpha^2(1-\alpha)+o_n(1)$ which is maximized when $\alpha=\frac{2}{3}$ to be $\frac{2}{9}+o_n(1)$.

\section{Almost all $0/1$ matrices are minimizers (sketch)}\label{sec:minimizers}

We sketch of a proof that a random $0/1$ $h \times h$ matrix is almost surely a  minimizer.
As mentioned in the introduction, the analogous case for graph inducibility (i.e., that almost all graphs are fractilizers) was proved in \cite{FHL-2017} and
an asymptotic version (that almost all graphs are almost fractilizers) was proved in \cite{yuster-2019}.
In both of these papers, the authors mention in the end of their papers that their methods can be extended to other combinatorial structures such as hypergraphs and tournaments after small number of relatively straightforward modifications. Indeed, we will sketch why this is also the case for $0/1$ matrices.
It should be emphasized that what follows is not a complete proof and the interested reader may look at the aforementioned papers making the appropriate (mostly straightforward) modifications to the matrix setting.
The key point, in both papers, is to prove that random graphs (in our case, random $0/1$ matrices) are far from having certain well-defined symmetries. These lack of symmetries limits the amount 
of overlap between two copies in the large graph (in our case, the large  matrix).
Let us next formalize some matrix notions equivalent to graph notions used in \cite{FHL-2017}.
Hereafter $H$ is an $h \times h$ $0/1$ matrix.

For two distinct rows $i,i'$ of $H$, let $s(i,i')$ denote the number of coordinates in which $i$ and $i'$ differ. We similarly define $s(j,j')$ for a pair of distinct columns. A subset $R$ of rows of $H$ is a {\em signature} if the $|R| \times h$ sub-matrix on these rows has all its columns distinct. Analogously, we can define signature for a subset of columns.
Call $H$ a {\em typical} matrix if the following conditions hold:
(i) For every subset $R$ of rows with $|R| \ge 0.99h$, for every subset $C$ of columns with $|C| \ge 0.99h$,
for any injective mapping $r$ of $R$ to another subset of rows with $R^* \subseteq R$ non-fixed points where $|R^*| \ge h/100$, and for any injective mapping $c$ of $C$ to another subset of columns, there are at least $|R^*||C|/100$ entries $(i,j)$ with $i \in R$ and $j \in C$ for which $H(i,j) \neq H(r(i),c(j))$.
(ii) analogously, for every subset $R$ of rows with $|R| \ge 0.99h$, for every subset $C$ of columns with $|C| \ge 0.99h$,
for any injective mapping $r$ of $R$ to another subset of rows, and for any injective mapping $c$ of $C$ to another subset of columns with $C^* \subseteq C$ non-fixed points where $|C^*| \ge h/100$, there are at least $|R||C^*|/100$ entries $(i,j)$ with $i \in R$ and $j \in C$ for which $H(i,j) \neq H(r(i),c(j))$.
(iii) For every pair of rows $i,i'$ it holds that $s(i,i') \ge h/4$ and similarly for every pair of columns $j,j'$ it holds that $s(j,j') \ge h/4$.

It is fairly standard to prove that as $h$ goes to infinity, a random $h \times h$ $0/1$ matrix $H$ (namely, where each entry of $H$ is independently chosen to be $1$ or $0$ equally likely), is typical with probability approaching $1$.

The following can then be obtained adapting the theorem from \cite{FHL-2017} to matrices:
For all $h$ sufficiently large, every $h \times h$ $0/1$ typical matrix is a minimizer.
A very useful fact in the proof is that in a typical matrix, there are small row/column signatures of order
only $O(\log h)$ (this easily follows from condition (c) in the definition of typical).
The key argument proceeds as follows: Suppose $M$ is an $n \times n$ matrix maximizing the number
of copies of $H$. Partition the row index set $R$ to $R_1,\ldots,R_h$ as follows.
Let $R_k$ be those rows for which $v \in R_k$ if the number of copies of $H$ in which row $v$ corresponds to row $k$ of $H$ is the largest over all choices of $k$ (note that this makes $R_1,\ldots,R_k$ a covering; make it a partition by arbitrarily breaking ties). Analogously partition the the column index set $C$ to
$C_1,\ldots,C_h$. An entry $(u,v) \in R_i \times C_j$ of $M$ is {\em consistent} if $M(u,v)=H(i,j)$, otherwise it is inconsistent. The key step is then to show that each inconsistent pair can only appear in a small number of copies of $H$ in $M$. Once this is established, it can be shown that if we flip the entries of inconsistent pairs we will obtain
a matrix with more $H$-copies, contradicting $M$'s maximality. This, in turn, gives that $M$ is
a blowup of $H$ with all entries in the sub-matrix $R_i \times C_j$ the same and equal $H(i,j)$.
A simple argument shows that the number of copies of $H$ in such a case is maximized when all part orders are $\lfloor n/h \rfloor$ or $\lceil n/h \rceil$.
To show that each inconsistent pair can only appear in a small number of copies of $H$ in $M$ (this is the most technical part) one uses the fact that there are small signatures and the typicality of $H$.

\section*{Acknowledgment}

The author thanks Asaf Shapira for useful discussions and Nati Linial for suggesting the use of
variational calculus.

\bibliographystyle{abbrv}

\bibliography{references}

\newpage

\appendix

\section{Kenyon's proof for \eqref{e:func3}}\label{sec:var}

For completeness we describe, almost verbatim (but more verbose for the sake of readers less familiar with variational calculus over probability measures) Kenyon's solution to \eqref{e:func3}.
Recall that we consider the functional 
\begin{equation}
	{\mathcal H}(h^*) = \int_{0}^{1} \int_{x}^{1} x(y-x)h^*(x)h^*(y) \, dy \, dx
\end{equation}
and wish to maximize it over all probability measures $h^*$ in $[0,1]$.
In particular, $h^* \ge 0$ in this range and $\int_0^1 h^*(x)dx = 1$.
Also take note that, being an arbitrary probability measure, $h^*$ is a generalized function.
Namley, at a point $t$ where $h^*$ has nonzero mass $p$, the value of $h^*$ is $p \delta_t$ where $\delta_t$ is the delta function at $t$ (so integrating over $t$ adds $p$ to the total).

As usual in variation calculus, we consider small perturbations of $h^*$ of the form $\varepsilon \eta(x)$. As we are replacing $h^*$ with its perturbation $h^*+\varepsilon \eta$, we are confined
to the initial condition that the total integral stays $1$, i.e., $\eta$ must be integrable and such that $\int_{0}^{1}\eta(x)dx=0$. We also need $\eta$ to be nonnegative at points where $h^*$ vanishes.

Let  us consider ${\mathcal H}(h^*+\varepsilon \eta)$ as a function of $\varepsilon$ and notice that it is differentiable. Since for a maximizing $h^*$ the value of the derivative at $\varepsilon = 0$ is $0$, we obtain
\begin{equation}\label{e:eps=0}
0 = \frac{d{\mathcal H}(h^*+\varepsilon \eta)}{d \varepsilon} \bigg\rvert_{\varepsilon = 0} =
\int_{0}^{1} \int_{x}^{1} x(y-x)(h^*(x)\eta(y)+\eta(x)h^*(y)) \, dy \, dx\;.
\end{equation}
As the latter holds for any $\eta$ satisfying the initial conditions, it in particular holds
for the following. Let $a$ be a point where $h^*(a) > 0$. Take $\eta=\delta'_a$ where we recall that the latter is formally defined through its integral property as
\begin{equation}\label{e:dfd}
\int_{-\infty}^{-\infty} f(u)\delta'_a(u)du = -f'(a)
\end{equation}
for any differentiable function $f$ on the entire real line (sometimes $f$ is referred to as a {\em test function}). In variational calculus,
$\delta'_a(u)$ is referred to as the {\em delta function derivative} and \eqref{e:dfd} is easily obtained via integration by parts from the fundamental property of the delta function
$\int_{-\infty}^{\infty}f(u)\delta_a(u)du = f(a)$.
Using $\eta=\delta'_a$ in each part of the r.h.s. of \eqref{e:eps=0} we obtain:
\begin{align*}
 & \int_{0}^{1} \int_{x}^{1} xyh^*(x)\delta'_a(y) \, dy \, dx\\
 = & \int_{0}^{1} xh^*(x) \left[\int_{x}^{1} y\delta'_a(y) \, dy\right] \, dx\\
 = &  \int_{0}^{1} xh^*(x) \left[-{\bf 1}_{[0,a]}(x)\right] \, dx\\
 = & - \int_{0}^{a} xh^*(x) \, dx\;.
\end{align*}
Similarly, we obtain
\begin{align*}
	\int_{0}^{1} \int_{x}^{1} -x^2h^*(x)\delta'_a(y) \, dy \, dx & = 0\,,\\	
	\int_{0}^{1} \int_{x}^{1} xy h^*(y)\delta'_a(x) \, dy \, dx & = \int_{a}^{1}-yh^*(y)\, dy\,,\\	
	\int_{0}^{1} \int_{x}^{1} -x^2 h^*(y)\delta'_a(x) \, dy \, dx & = \int_{a}^{1}2a h^*(y)\, dy\,.	
\end{align*}
So, altogether we obtain from \eqref{e:eps=0} that
$$
 0 = - \int_{0}^a xh^*(x)\,dx + \int_{a}^1 (2a-y)h^*(y)\,dy\;.
$$
Now, using integration by parts and letting $H^*(a) = \int_{0}^a h^*(x)\,dx$ denote the cumulative distribution function of $h^*$, we obtain from the last equality that
$$
\int_{0}^1 xh^*(x)\,dx = 2a(H^*(1)-H^*(a))\;.
$$
Denoting the left hand side (which is a constant depending only on $h^*$) by $2c$, we obtain that
$H^*(a)=H^*(1)-c/a$ and hence it must be that $h^*(a)=c/a^2$ whenever $h^*(a) > 0$.
Notice however, that our argument breaks if it happens to hold that $h^*(a) > 0$ yet $h^*$ vanishes on
an interval contained in $(a,1]$. We shall prove that this cannot happen. So, assume for contradiction that $b > a$ and $b$ is contained in an interval where $h^*$ vanishes.
Let $\eta = \delta_b - \delta_a$ and observe that this is a valid perturbation since the integral
of such a choice of $\eta$ is still $1-1=0$ and is not negative in points where $h^*$ vanishes.
Considering the parts of the r.h.s. of \eqref{e:eps=0} for this choice we obtain
\begin{align*}
	\int_{0}^{1} \int_{x}^{1} xyh^*(x)(\delta_b(y)-\delta_a(y)) \, dy \, dx & =
	(b-a) \int_{0}^{a} xh^*(x)dx + b \int_{a}^{b} xh^*(x)dx\,.\\
\int_{0}^{1} \int_{x}^{1} -x^2h^*(x)(\delta_b(y)-\delta_a(y)) \, dy \, dx & =
 \int_{a}^{b} -x^2h^*(x)dx\,.\\
\int_{0}^{1} \int_{x}^{1} xyh^*(y)(\delta_b(x)-\delta_a(x)) \, dy \, dx & =
-a\int_{a}^{b} y h^*(y) dy + (b-a)\int_{b}^{1}y h^*(y) dy\,.\\
\int_{0}^{1} \int_{x}^{1} -x^2h^*(y)(\delta_b(x)-\delta_a(x)) \, dy \, dx & =
a^2\int_{a}^{b} h^*(y) dy + (a^2-b^2)\int_{b}^{1} h^*(y) dy\,.
\end{align*}
Summing all four parts we obtain from \eqref{e:eps=0} that
$$
0 = (b-a)\int_{0}^{1}xh^*(x)dx + \int_{a}^b (-x^2+a^2)h^*(x) dx +(a^2-b^2)\int_{b}^{1}h^*(x)dx\;.
$$
Considering the r.h.s. of the last equality as a function of $b$, we have that the middle integral in independent of $b$ near $b$ since $h^*$ vanishes around $b$, but then we have
$0 = 2c(b-a) + D + (a^2-b^2)(H^*(1)-H^*(b))$ for some $D$ which is a constant  near $b$,
implying that $G(b)$ is nonconstant around $b$, contradicting that $h^*$ vanishes around $b$.

So, we have that $h^*$ must be a probability measure in $[0,1]$ that is perhaps zero in an initial interval $[0,\alpha]$,
and then is of the form $c/x^2$ in $(\alpha,1)$. It may possibly have mass $\beta\delta_1$ at the point $1$. Summing up, we have
$$
	h^*(x) =
\begin{cases}
	0 				   & {\rm for~} 0 \le x \le \alpha\,, \\
	\frac{c}{x^2}	   & {\rm for~} \alpha < x < 1\,, \\
	\beta\delta_1(x)   & {\rm for~}  x = 1 \,.
\end{cases}	
$$
Using the initial condition $\int_{0}^1 h^*(x)dx = 1$ and the condition $2c = \int_{0}^1 xh^*(x) dx$
we can easily optimize for maximizing ${\mathcal H}(h^*)$ and obtain $\alpha=\beta=c=1/e$
for which we obtain ${\mathcal H}(h^*)=1/2e^2$, as claimed.

\section{Upper bound for $f(\B)$}\label{sec:flag}

Our upper bounds for $f(\B)$ is via flag algebra.
This method, pioneered in a classical paper of Razborov \cite{razborov-2007}, has become a
widely used indispensable tool in the area of homomorphism density problems in extremal combinatorics and other related areas; see \cite{razborov-2013} for a survey of flag algebra applications. To use this method in our setting, we first need to introduce the flag algebra notation and objects that we work with, which we describe in some brevity (yet in full).
The reader interested in more details may consult various surveys and gentle treatments to the subject, such a the one in \cite{GGHLM-2022}.

Throughout the remainder of this section we assume that all matrices are $0/1$, but note that all arguments can easily be generalized to other symbol sets.
We start with a simple observation which follows directly by double counting.
Let $\Set_{r,c}$ be the set of
$r \times c$ matrices. Suppose that $h \le r \le n$ and $h \le c \le n$.
By double counting, we have for $H^* \in \Set_{h,h}$ and $M \in \Set_{n,n}$ that
\begin{equation}\label{e:dc}
	d(H^*,M) = \sum_{H \in \Set_{r,c}}d(H^*,H)d(H,M)\;.
\end{equation} 

A {\em type} $\sigma$ is just an $a \times b$ matrix. Let $[\sigma]=(a,b)$ be the {\em size pair} of $\sigma$.
For a type $\sigma$ having size pair $(a,b)$, a {\em $\sigma$-flag} $F$ is a triple $(M, \theta_1,\theta_2)$ where $M$ is 
an $r \times c$ matrix, $\theta_1 \in \binom{[r]}{a}$
and $\theta_2 \in \binom{[c]}{b}$. We think of $\theta_1$ and $\theta_2$ as marking some of the rows and some of the columns, respectively. Furthermore, we require that the submatrix of $M$ corresponding to the marked rows
and columns is equal to $\sigma$.
Denote by $\Flag_{r,c}^\sigma$ the set of $r \times c$\,  $\sigma$-flags whose underlying matrix
is in $\Set_{r,c}$.
Figure \ref{f:flags} exhibits some $\sigma$-flags where $\sigma=(0,1)$ (hence $[\sigma]=(1,2)$).

\begin{figure}[h]
	\[
	\begin{blockarray}{cccc}
		& \bullet & \bullet & \\
		\begin{block}{p{1pt}(ccc)}
			$\bullet$ & 0 & 1 & 0\\
			& 1 & 0 & 1\\
		\end{block}
	\end{blockarray}
	\quad
	\begin{blockarray}{cccc}
		& \bullet & \bullet & \\
		\begin{block}{p{1pt}(ccc)}
			$\bullet$ & 0 & 1 & 1\\
			& 1 & 0 & 1\\
		\end{block}
	\end{blockarray}
	\quad
	\begin{blockarray}{cccc}
		& \bullet & & \bullet \\
		\begin{block}{p{1pt}(ccc)}
			& 0 & 1 & 0\\
			$\bullet$ & 0 & 0 & 1\\
		\end{block}
	\end{blockarray}
	\quad
	\begin{blockarray}{cccc}
		& & \bullet & \bullet \\
		\begin{block}{p{1pt}(ccc)}
			$\bullet$ & 0 & 0 & 1\\
			& 1 & 0 & 1\\
		\end{block}
	\end{blockarray}
	\]
	\caption{Some elements of $\Flag_{2,3}^\sigma$ where $\sigma=(0,1)$.}\label{f:flags}
\end{figure}

A {\em subflag} of a $\sigma$-flag $(M,\theta_1,\theta_2)$ is a $\sigma$-flag corresponding to a submatrix of $M$ whose set of rows contains
all marked rows (and possibly some other rows) and whose set of columns contains all marked columns (and possibly some other columns).
Given $\sigma$-flags $F \in \Flag_{r,c}^\sigma$ and $K \in \Flag_{m,n}^\sigma$, let
$p(F,K)$ be the probability that a randomly chosen $r \times c$ subflag of $K$ equals $F$
(if $m < r$ or $n < c$ or $K$ is not a $\sigma$-flag, define $p(F,K)=0$). 

Similarly, given flags $F,F' \in \Flag_{r,c}^\sigma$ and $K \in \Flag_{m,n}^\sigma$,
define the {\em joint density} $p(F,F';K)$ as the probability that if we randomly choose a
pair of $r \times c$ subflags $(S,S')$ of $K$ subject to the set of unmarked rows of
$S$ being disjoint from the unmarked rows of $S'$ and the unmarked columns of
$S$ being disjoint from the unmarked columns of $S'$, then $(S,S')=(F,F')$.
If $K$ is not a $\sigma$-flag or $m < 2r-a$, or $n < 2c-b$ where $[\sigma]=(a,b)$, then define $p(F,F';K)=0$. As an illustrative example of these notions, consider the type $\sigma$, and the $\sigma$-flags $F,F',K$ with the corresponding values of $p(F,K)$ and $p(F,F';K)$ given in Figure \ref{f:densities} and its caption.

\begin{figure}[ht]
	\[
	\sigma=(0)
	\qquad
	F = 
	\begin{blockarray}{ccc}
		& \bullet & ~ \\
		\begin{block}{p{1pt}(cc)}
			& 1 & 0\\
			$\bullet$ & 0 & 1\\
		\end{block}
	\end{blockarray}
	\qquad
	F' = 
	\begin{blockarray}{ccc}
		&  & \bullet \\
		\begin{block}{p{1pt}(cc)}
			$\bullet$ & 1 & 0\\
			& 1 & 1\\
		\end{block}
	\end{blockarray}
	\qquad
	K = 
	\begin{blockarray}{ccccc}
		&  & \bullet & & \\
		\begin{block}{p{1pt}(cccc)}
			& 1 & 1 & 0 & 0\\
			$\bullet$ & 1 & 0 &1 & 1\\
			& 1 & 1 &0 & 0\\
		\end{block}
	\end{blockarray}
	\]
	\caption{Depicted are a $\sigma$-type and three $\sigma$-flags $F,F' \in \Flag_{2,2}^\sigma$,
		$K=(M,\{2\},\{2\}) \in \Flag_{3,4}^\sigma$. Notice that $K$ has $\binom{2}{1}\binom{3}{1}=6$\; $2 \times 2$ subflags.
		Only two of them equal $F$, so $p(F,K)=\frac{1}{3}$.
		Similarly, $K$ has $\binom{2}{1}\binom{3}{1}\binom{1}{1}\binom{2}{1}=12$ pairs of $2 \times 2$ subflags
		such that the unmarked rows in both are disjoint and the unmarked columns in both are disjoint.
		Only one such pair corresponds to $(F,F')$. The first, corresponding to $F$, is given by choosing rows $\{1,2\}$ and columns $\{2,3\}$ of $M$, and the second, corresponding to $F'$, is given by choosing rows
		$\{2,3\}$ and columns $\{1,2\}$ of $M$. Hence, $p(F,F';K)=\frac{1}{12}$.}\label{f:densities} 
\end{figure}

Now, suppose that $\Flag_{r^*,c^*}^\sigma = \{F_1,\ldots,F_t\}$.
Let $Q$ be a $t \times t$ positive semidefinite matrix and let $K \in \Flag^\sigma_{n,n}$.
It follows from positive-semidefiniteness and from Lemma 2.3 in \cite{razborov-2007} that:
$$
0 \le \sum_{i,j}Q[i,j]p(F_i,F_j;K)+o_n(1)\;.
$$ 
Notice that the last inequality holds also if $K$ is not a $\sigma$-flag as $p(F_i,F_j;K)=0$ in this case. 

Using notation analogous to \cite{BHL+-2015}, for a matrix $M$, let $\Theta((a,b),M)$ be the set of all pairs $(\theta_1,\theta_2)$
such that $\theta_1$ is a subset of $a$ row indices of $M$ and $\theta_2$ is a subset of
$b$ column indices of $M$. Suppose that $[\sigma]=(a,b)$ and consider a uniform distribution on
$\Theta((a,b),M)$. Then $K=(M,\theta_1,\theta_2)$ may or may not be a $\sigma$-flag.
In any case, we have from the last inequality that if $M$ is an $n \times n$ matrix, then: 
\begin{equation}\label{e:o1}
	0 \le \sum_{i,j}Q[i,j]{\mathbb E}_{(\theta_1,\theta_2) \in \Theta([\sigma],M)}[p(F_i,F_j;(M,\theta_1,\theta_2))]+o_n(1)\;.
\end{equation}

Recalling the definition of $\Set_{r,c}$ we have by double counting,
and assuming $r \ge 2r^*-a$ and $c \ge 2c^*-b$, that for $M \in \Set_{n,n}$ it holds that
$$
{\mathbb E}_{(\theta_1,\theta_2) \in \Theta([\sigma],M)}[p(F_i,F_j;(M,\theta_1,\theta_2))] =
\sum_{H \in \Set_{r,c}} {\mathbb E}_{(\theta_1,\theta_2) \in \Theta([\sigma],H)}[p(F_i,F_j;(H,\theta_1,\theta_2))] \cdot d(H,M)\;.
$$
Hence, together with \eqref{e:o1},
$$
0 \le \sum_{H \in \Set_{r,c}} \left(\sum_{i,j}Q[i,j] {\mathbb E}_{(\theta_1,\theta_2) \in \Theta([\sigma],H)}[p(F_i,F_j;(H,\theta_1,\theta_2))] \right)d(H,M) + o_n(1)\;. 
$$
Note that the coefficient of $d(H,M)$ is independent on $M$ as it only depends on
$\sigma,r^*,c^*,Q,H$, so denote it by $c_H(\sigma,r^*,c^*,Q)$. Rewriting, we have
$$
0 \le \sum_{H \in \Set_{r,c}}c_H(\sigma,r^*,c^*,Q)d(H,M)+o_n(1)\;.
$$ 

For $i \in [s]$, let $(\sigma_i,r_i,c_i,Q_i)$ be quadruples satisfying $r \ge 2r_i-a_i$
and $c \ge 2c_i-b_i$, where $[\sigma_i]=(a_i,b_i)$, and $Q_i$ is a positive semidefinite matrix
whose rows and columns are indexed by $\Flag_{r_i,c_i}^{\sigma_i}$.
For each quadruple, we consider the corresponding coefficient $c_H(\sigma_i,r_i,c_i,Q_i)$ and set
$c_H = \sum_{i=1}^s c_H(\sigma_i,r_i,c_i,Q_i)$. Summing the last displayed inequality for each quadruple we have:
$$
0 \le \sum_{H \in \Set_{r,c}}c_H \cdot d(H,M)+o_n(1)\;.
$$ 
The last inequality together with \eqref{e:dc} gives for $H^* \in \Set_{h,h}$ where
$h \le \min\{r,c\}$ that
$$
d(H^*,M) - o_n(1) \le \sum_{H \in \Set_{r,c}}(d(H^*,H)+c_H)d(H,M) \le \max_{H \in \Set_{r,c}}(d(H^*,H)+c_H)\;.
$$
We therefore obtain: 
\begin{corollary}\label{corr:flag}
	Let $h,s,r,c$ be positive integers such that $h \le \min\{r,c\}$.
	Let $H^* \in \Set_{h,h}$.
	Let $(\sigma_i,r_i,c_i,Q_i)$ for $i \in [s]$
	be such that $\sigma_i$ is a type where $[\sigma_i]=(a_i,b_i)$ such that $r_i > a_i$ and $c_i > b_i$ satisfy $r \ge 2r_i-a_i$, $c \ge 2c_i-b_i$ and $Q_i$ is a positive semidefinite matrix indexed by $\Flag_{r_i,c_i}^{\sigma_i}$. Then,
	\[
	\pushQED{\qed}
	f(H^*) \le \max_{H \in \Set_{r,c}}(d(H^*,H)+c_H)\;. \qedhere \popQED
	\] 
\end{corollary}
By Corollary \ref{corr:flag}, to obtain an upper bound for $f(\B)$
we can choose integers $r,c \ge 2$,
compute $d(\B,H)$ for all $H \in \Set_{r,c}$, choose $s$, choose types $\sigma_i$ and sizes $r_i,c_i$ for
$i \in [s]$ and then solve a semidefinite program to
obtain solutions for the semidefinite matrices $Q_i$ for $i \in [s]$ that minimizes
$\max_{H \in \Set_{r,c}}(d(\B,H)+c_H)$.
This is detailed in the following subsection. 

\subsection{Choosing types and flags}\label{subsec:ub}

The results in this subsection reference a computer program which we refer to as the ``generator program'' (link to code given in Table {\ref{table:urls}). We refer to the notation established in Section \ref{sec:flag} above.
	
	We shall use $H^*=\B$ and use $r=3$ and $s=4$.
	The generator program constructs all elements in $\Set_{3,4}$.
	Note that, trivially, $|\Set_{3,4}|=4096$ as it contains all possible binary $3 \times 4$ matrices. The generator program next computes $d(\B,H)$ for every $H \in \Set_{3,4}$.
	
	We shall use the $1 \times 2$ types $\sigma_1=(0,0)$, $\sigma_2=(0,1)$, $\sigma_3=(1,1)$ and $\sigma_4=(1,0)$, so we have $s=4$.
	We shall use the flag list $\Flag_{2,3}^{\sigma_i}$ for each of these types.	
	The generator program computes the elements in these flag lists.
	We have that $|\Flag_{2,3}^{\sigma_i}|=96$ for each of the four types. This is easy to see since the type forces choosing two cells
	in the same row (six choices) and this leaves four cells that can be occupied freely with $0$ or $1$.
	
	Finally, the generator program computes the last piece of data needed to generate the semidefinite program, namely, for each type $\sigma_i$, for each pair of flags $F,F' \in \Flag_{2,3}^{\sigma_i}$, and for each
	$H \in \Set_{3,4}$ it computes
	${\mathbb E}_{(\theta_1,\theta_2) \in \Theta([\sigma_i],H)}[p(F,F';(H,\theta_1,\theta_2))]$.
	
	Using the computed constants mentioned above, the generator program then creates an input file in standard format (see next subsection for more details) that is given as input to an sdp solver.
	The solver gives corresponding positive semidefinite matrices $Q_i$ for $i \in [s]$, for which we
	obtain that
	$$
	f(\B)  \le \max_{H \in \Set(\emptyset)_{3,4}}(d(\B,H)+c_H) \approx 0.2716956... \,.
	$$ 
	All of this bounds are proved rigorous using a rounding procedure yielding a rational certificate elaborated upon Subsection \ref{subsec:rounding}.
	
	\subsection{Converting to standard sdp format}
	
	State of the art numerical sdp solvers such as CSDP \cite{borchers-1999} (which is the one we use) inevitably use floating point arithmetic and hence their results, while highly accurate, are, nevertheless, approximations. One then needs to apply some further rounding procedure
	in order to turn the results into a rigorous proof. See Subsection \ref{subsec:rounding} for details on our rounding method of choice.
	
	CSDP solves the following standard form semidefinite program.
	Let $C,A_1,\ldots,A_m$ be given real symmetric matrices and let $X$ be a (variable) real symmetric matrix.
	Let $a=(a_1,\ldots,a_m)$ be a given real vector.
	The semidefinite program solved by CSDP is:
	\begin{equation}\label{e:csdp}
		\begin{aligned}
			\max \hspace{32pt} & \operatorname{tr}(C X) \\
			\text{subject to } & \operatorname{tr}(A_jX) = a_j & \text{for } 1 \le j \le m\;,\\
			& X \succeq 0
		\end{aligned}
	\end{equation}
	(here $X \succeq 0$ means $X$ is positive semidefinite).
	
	Recalling the notation in Section \ref{sec:flag} and Corollary \ref{corr:flag}, note that our  semidefinite program is:
	$$
	\min_{Q_1,\ldots,Q_s} \max_{H \in \Set_{3,4}}(d(\B,H)+c_H)
	$$
	which is equivalent to the program
	$$
	1 - \max_{Q_1,\ldots,Q_s} \min_{H \in \Set_{3,4}}(1-d(\B,H)-c_H)
	$$
	subject to $Q_i \succeq 0$ and to $c_H = \sum_{i=1}^s c_H(\sigma_i,2,3,Q_i)$.
	Denoting ${\mathcal F}_{2,3}^{\sigma_i}=\{F^i_1,\ldots,F^i_{t_i}\}$ recall 
	the definition of $c_H(\sigma_i,2,3,Q_i)$ given in Section \ref{sec:flag}:
	$$
	c_H(\sigma_i,2,3,Q_i) = \sum_{u=1}^{t_i}\sum_{v=1}^{t_i}
	Q_i[u,v]{\mathbb E}_{(\theta_1,\theta_2) \in \Theta([\sigma_i],H)}[p(F^i_u,F^i_v;(H,\theta_1,\theta_2))]\;.
	$$
	Our aim is therefore to solve the program
	$$
	\max_{Q_1,\ldots,Q_s} \min_{H \in \Set_{3,4}}(1-d(\B,H)-c_H)
	$$
	subject to $Q_i \succeq 0$. Now, observe that if $M$ and $N$ are positive integers, then
	the last program yields the same result as the program
	$$
	\frac{1}{M}\max_{Q_1,\ldots,Q_s} \min_{H \in \Set_{3,4}}(M \cdot(1-d(\B,H))-Nc_H)\;.
	$$
	The advantage of using the latter formulation is that we can choose $M$, $N$ such that all constants in the
	sdp are integers, which makes the input to the sdp more concise.
	Since $d(\B,H)$ is a density of a $2 \times 2$ matrix inside at $3 \times 4$ matrix, it is a rational with denominator $18$, so we shall choose $M=18$.
	Also notice that for $i \in [s]$, $p(F^i_u,F^i_v;(H,\theta_1,\theta_2))]$ is a rational with denominator
	$72$ (the number of choices for $(\theta_1,\theta_2) \in \Theta([\sigma_i],H)$ is $\binom{3}{1}\binom{4}{2}=18$ and the number of choices for a subset pair $(S,S')$ in the definition of joint density is $2 \times 2 = 4$ in our case), so we choose $N=72$.
	We then execute the sdp and obtain by Corollary \ref{corr:flag} that
	\begin{equation}\label{e:sdp}
		f(\B)) \le 1 - \frac{1}{18}\max_{Q_1,\ldots,Q_s} \min_{H \in \Set_{3,4}}(18(1-d(\B,H))-72c_H)\;.
	\end{equation}
	
	Translating the sdp \eqref{e:sdp} to the standard sdp format of \eqref{e:csdp} is a straightforward process of adding slack variables. 
	In our case $m=|\Set_{3,4}|$ is the dimension of the vector $a=(a_1,\ldots,a_m)$.
	Now, suppose that $\Set_{3,4}=\{H_1,\ldots,H_m\}$, then $a_i=18(1-d(\B,H_i))$.
	Each of the matrices $C,A_1,\ldots,A_m$ is a block diagonal matrix with precisely $s+1$ blocks, where
	for $i \in [s]$, the $i$'th block corresponds to the type $\sigma_i$ and the last block is the ``slack block''.
	For $i \in [s]$, the order of block $i$ in each of these matrices is $|{\mathcal F}_{2,3}^{\sigma_i}|$
	(namely, it is $96$).
	The dimension of the slack block is $m+1$.
	As for their entries, the matrix $C$ is entirely zero except for the $[1,1]$ entry of the slack block which is $1$.
	For $1 \le i \le [s]$ and for $1 \le j \le m$, entry $[u,v]$ of block $i$ of $A_j$ is
	$72{\mathbb E}_{(\theta_1,\theta_2) \in \Theta([\sigma_i],H_j)}[p(F^i_u,F^i_v;(H_j,\theta_1,\theta_2))]$.
	The slack block of $A_j$ is entirely zero except for entries $[1,1]$ and $[j+1,j+1]$ which are $1$
	(so the slack block is a diagonal matrix).
	The input files referenced in Table \ref{table:urls} contain all of these values in standard
	SDPA sparse format (see the manuals of either SDPA \cite{YFK-2003} or CSDP for a description of this format).
	
	\subsection{Rounding}\label{subsec:rounding}
	
	As the output matrix of the sdp solver is an approximate floating point solution, one needs to couple the approximate result with a rounding argument in order to obtain a rigorous proof. We will adopt a rounding method similar to the one in \cite{SS-2018}.
	
	Let $X'$ denote the result matrix of the sdp execution.
	Notice that $X'$ has the same block structure as the matrices $A_j$. Recalling that
	block $i$ is square symmetric of order $|{\mathcal F}_{2,3}^{\sigma_i}|$ for $i \in [s]$
	and that the slack block is diagonal of order $m+1=|\Set_{3,4}| + 1$, we have that the number of variables in our problem is:
	$$
	\sum_{i=1}^s \binom{|{\mathcal F}_{2,3}^{\sigma_i}|+1}{2} + |\Set_{3,4}| + 1 = 22721\;.
	$$
	
	Let $\delta > 0$ be a parameter, which can be configured by CSDP and is taken by default to be $10^{-8}$
	\cite{borchers-1999}. CSDP guarantees \cite{borchers-1999} that $X'$ satisfies:
	\begin{equation}\label{e:guarantee}
		\begin{aligned}
			|\operatorname{tr}(A_jX') - a_j| & \le \delta & \text{for } 1 \le j \le m \;,\\
			X' & \succeq 0\;.
		\end{aligned}
	\end{equation}
	Nevertheless, as $X'$ consists of floating point numbers, we need to convert them to rationals in order to certify a rigorous proof, without incurring significant loss in the obtained bound.
	Since $X' \succeq 0$, we use python numpy (which uses BLAS/LAPACK \cite{anderson-1999})  to compute the
	Cholesky decomposition of $X'$, so $L'L'^T=X'$ where $L'$ is lower triangular. Nevertheless, the computed $L'$ is still a floating point approximation (since $X'$ is an approximation and since the Cholesky decomposition routine may further incur additional loss). Let $D$ be a large integer (specifically, we use $D=10^6$).
	We multiply each entry of $L'$ by $D$ and round each resulting entry to the closest integer, thus we obtain
	an integer matrix $L$. We expect that $\frac{1}{D^2}LL^T$ is a good approximation of the exact result $X$
	of our sdp. The matrix $L$ computed by our python script (see Table \ref{table:urls}) is our certificate and is provided by the links in Table \ref{table:urls} as well.
	
	Our final task is to verify that $L$ produces a bound close to optimal. Let $M=LL^T$. We first verify that $L$ is indeed the Cholesky decomposition of $M$. Since $L$ is lower triangular, this means that we just need
	to verify that all diagonal entries of $L$ are positive. Indeed, as our script shows, this holds
	(the lowest diagonal entry of $L$ is $6$).
	
	Let $b_j = \operatorname{tr}(\frac{1}{D^2}A_jM-a_j)$. Since $A_j,M,a_j$ are all integral, we have that $b_j$
	is a rational with denominator $D^2$. Let $\varepsilon = \max_{j=1}^m|b_j|$. Our program shows that
	the obtained {\em rational} $\varepsilon$ is smaller than $3.5\cdot 10^{-5}$.
	
	Notice that $\operatorname{tr}(\frac{1}{D^2}CM)$ is a lower bound for the {\em precise solution} of the problem
	\begin{equation}\label{e:csdp-precise}
		\begin{aligned}
			\max \hspace{32pt} & \operatorname{tr}(CX) \\
			\text{subject to } & \operatorname{tr}(A_jX) = a_j+b_j & \text{for } 1 \le j \le m\;,\\
			& X \succeq 0\;.
		\end{aligned}
	\end{equation}
	Recalling that $a_j=18(1-d(\B,H_i))$ we have, in turn, that the solution of \eqref{e:csdp-precise} is the solution of 
	$$
	\max_{Q_1,\ldots,Q_s} \min_{H_j \in \Set_{3,4}}(18(1-d(\B,H_i))+b_j-72c_H)\;.
	$$
	Hence, if $Q_1,\ldots,Q_s$ denote the non-slack blocks of $\frac{1}{D^2}M$, we have that
	$$
	\operatorname{tr}\left(\frac{1}{D^2}CM\right) \le \min_{H_j \in \Set_{3,4}}(18(1-d(\B,H_i))+b_j-72c_H)\;.
	$$
	Since $|b_j| \le \varepsilon$, we have that
	$$
	\operatorname{tr}\left(\frac{1}{D^2}CM\right) - \varepsilon \le \min_{H_j \in \Set_{3,4}}(18(1-d(\B,H_i))-c_H)\;.
	$$
	Hence, by \eqref{e:sdp}, $18-\operatorname{tr}(\frac{1}{D^2}CM)+\varepsilon \ge 18f(\B)$.
	Finally, running our python script referenced in Table \ref{table:urls}, or just examining entry
	$[1,1]$ of the slack block of $M$ (recall that the only nonzero entry of $C$ is the $[1,1]$ entry of the slack block, which equals $1$), we have
	$$
	\operatorname{tr}(CM) = 13109511938436 \;.
	$$
	Using the aforementioned value of $\varepsilon$ and that $D=10^6$, we have
	$f(\B) \le 0.2716956...$ so the bounds claimed in Subsection \ref{subsec:ub} hold.
	\begin{table}[ht]
		\centering
		\begin{tabular}{l|l}
			\hline
			Generator program & \footnotesize{\url{github.com/raphaelyuster/matrix-density/blob/main/matrix-density.cpp}} \\
			SDP rounding script & \footnotesize{\url{github.com/raphaelyuster/matrix-density/blob/main/sdp_analyzer.py}} \\
			SDP input & \footnotesize{\url{raw.githubusercontent.com/raphaelyuster/matrix-density/main/emptyset.dat-s}} \\
			Certificate & \footnotesize{\url{github.com/raphaelyuster/matrix-density/blob/main/emptyset.cert}} \\
			Script output& \footnotesize{\url{github.com/raphaelyuster/matrix-density/blob/main/emptyset.txt}} \\
			\hline
		\end{tabular}
		\caption{Description of each certificate and its corresponding url.}
		\label{table:urls} 
	\end{table}  

\end{document}